\documentclass[]{amsart}
\usepackage{amsmath}
\usepackage{amsfonts}
\usepackage{MnSymbol,mathrsfs}
\usepackage{stmaryrd}

\usepackage{amscd}
\usepackage[all]{xy}

\theoremstyle{plain}
\newtheorem{prop}[]{Proposition}
\newtheorem{lem}[prop]{Lemma}
\newtheorem{thm}[]{Theorem}
\newtheorem{cor}[prop]{Corollary}
\theoremstyle{definition}
\newtheorem{dfn}[prop]{Definition}
\theoremstyle{remark}
\newtheorem{rmk}[prop]{Remark}
\newtheorem{example}[]{Example}

\newcommand{\ensemble}[1]{\left\{ #1 \right\}}
\newcommand{\suchthat}{\mid}
\newcommand{\norm}[1]{\lVert #1 \rVert}
\newcommand{\bnorm}[1]{{\left\vert\kern-0.25ex\left\vert\kern-0.25ex\left\vert #1 \right\vert\kern-0.25ex\right\vert\kern-0.25ex\right\vert}}
\newcommand{\banach}[1]{\mathrm{ #1 }}
\newcommand{\absolute}[1]{\left| #1 \right|}
\newcommand{\pairing}[2]{\left \langle #1, #2 \right\rangle}

\newcommand{\schwfuncs}{\mathcal{S}}
\newcommand{\N}{\mathbb{N}}
\newcommand{\Z}{\mathbb{Z}}
\newcommand{\Q}{\mathbb{Q}}
\newcommand{\C}{\mathbb{C}}
\newcommand{\R}{\mathbb{R}}
\newcommand{\T}{\mathbb{T}}
\newcommand{\Tang}{\mathrm{T}}

\newcommand{\smooth}[1]{\mathcal{#1}}
\newcommand{\Mult}{\mathcal{M}}

\DeclareMathOperator{\ind}{Ind}

\DeclareMathOperator{\ch}{ch}
\DeclareMathOperator{\Todd}{Todd}

\DeclareMathOperator{\Id}{Id}
\DeclareMathOperator{\ptensor}{\hat{\otimes}}

\DeclareMathOperator{\HP}{HP}
\DeclareMathOperator{\HC}{HC}
\newcommand{\KK}{\ensuremath{\mathit{KK}}}
\DeclareMathOperator{\Det}{Det^{+}_\tau}
\DeclareMathOperator{\Met}{Met_\tau}
\DeclareMathOperator{\SMet}{SMet^{(\omega)}_\tau}
\DeclareMathOperator{\topol}{top}
\DeclareMathOperator{\SL}{SL}
\DeclareMathOperator{\Sp}{Sp}
\DeclareMathOperator{\GL}{GL}

\usepackage{hyperref}

\usepackage[shortalphabetic,initials,lite,msc-links,nobysame]{amsrefs}

\newcommand{\condprojinv}{\ensuremath{\mathsf{(PI)}}}
\newcommand{\conduniterg}{\ensuremath{\mathsf{(UE)}}}

\begin{document}

\title[Projectively Measured Foliation]{Operator algebra of foliations with projectively invariant transverse measure}
\author[M. Yamashita]{Makoto Yamashita}
\address{Graduate School of Mathematical Sciences, the University of Tokyo}
\curraddr{Department of Mathematics, Ochanomizu University}
\email{makotoy@ms.u-tokyo.ac.jp}

\begin{abstract}
  We study the structure of operator algebras associated with the foliations which have projectively
  invariant measures.  When a certain ergodicity condition on the measure
  preserving holonomies holds, the lack of holonomy
  invariant transverse measure can be established in terms of a cyclic cohomology
  class associated with the transverse fundamental cocycle and the modular
  automorphism group.
\end{abstract}

\maketitle

\section{Introduction}
\label{sec:intro}

In~\cite{MR636521}, A. Connes introduced the von Neumann algebra $W(M; F)$
associated with a foliation $F$ on a manifold $M$.  He showed that several measure theoretic properties of the space of leaves $M/F$ can be stated in terms of $W(M; F)$.  For example, the ergodic components of $M/F$ correspond to the continuous decomposition of $W(M; F)$, and the lack of holonomy invariant
measure in the absolute continuity class of transverse Lebesgue measure is equivalent to that $W(M; F)$ being a type III algebra.  Subsequently it was shown by
S. Hurder and A. Katok~\cite{MR752795}, and Connes~\cite{MR866491} that if a foliation $F$ admits a nonzero generalized Godbillon--Vey class $\omega$, then the direct summand of $W(M; F)$ corresponding to the support of $\omega$ is of type III.

When one considers the holonomy of transverse coordinates, there is an interesting class of foliations with simplified transverse
structure, namely that of transversely affine foliations.  For such ones,
several sufficient conditions for the existence of an invariant transverse
measure were established by W. M. Goldman, M. W. Hirsch, and G. Levitt~\cite{MR671227}.  There are
also studies of more specific examples in this class, for example the ones given
by the hyperbolic automorphisms of tori, by J. F. Plante~\cite{MR609116}.

In this paper we investigate the class foliation with positive projectively invariant
measures, naturally including the transversely affine foliations.  Our main theorem gives a sufficient condition (Theorem~\ref{thm:inv-meas-flow-weights}) for $W(M; F)$ to be of type
III in terms of a certain cyclic cocycle $i_D\phi$ and its pairing with the $K$-group of the foliation algebra.  Although we closely follow Connes's method to establish the lack of invariant transverse measure, we note that the Godbillon--Vey classes themselves can be trivial in our setting (Remark~\ref{rmk:godbillon-vey-num-triv}).
We also obtain a more precise description (Theorem
\ref{thm:iii-lambda-proj-inv-per-calc-trans-fund}) of the pairing of $i_D\phi$
with $K_*(C^*_r(M; F))$ when $W(M; F)$ is of type
III$_\lambda$ for $0 < \lambda < 1$.

This paper is organised as follows.  Section
\ref{sec:preliminaries} is devoted to a review of basic groupoid and algebra
constructions associated with foliation.  Next in Section~\ref{sec:peri-acti-cycl}, we give a description
(Corollary~\ref{cor:period-dual-action-pairing-der-cycle-k-group}) of the image of
invariant cyclic cocycles under the `boundary map' in the Pimsner--Voiculescu
type exact sequence of periodic cyclic cohomology by R. Nest~\cite{MR961899}.  In Section~\ref{sec:invar-fund-class}, we show that the transverse fundamental
cocycle $\phi$ is invariant under the modular automorphism group if the
restricted holonomy groupoid admits a projective invariant density.  Due to this
invariance, one obtains a cyclic cocycle of degree $q + 1$ over the
smooth convolution algebra of $F$ as the `interior product' $i_D \phi$ of the
generator of the modular automorphism group with $\phi$.  This cocycle is anabelian in the sense of Connes, and its pairing with a class in the $K$-group gives rise to a invariant measure on the flow of weights, leading to the theorems mentioned above.

\section{Preliminaries}
\label{sec:preliminaries}

\subsection{Basic definitions of foliation  algebras}

Let $M$ be a smooth oriented manifold of dimension $n$, and $F$ a smooth
oriented foliation of dimension $p$ on $M$.  We let $\Tang M$ denote the tangent bundle of $M$ and identify $F$ with an integrable
subbundle of $\Tang M$.  The normal bundle $\Tang M / F$ of $F$ is denoted by $\tau$, and the codimension $n - p$ of
$F$ is by $q$.  Thus $\wedge^q \tau^*$ can be regarded as the bundle of signed transverse densities.  Let $\Det(M)$ denote the bundle of nonzero transverse densities over $M$, which is a subbundle of $\wedge^q \tau^*$ determined by the orientation on $F$.  It is a principal $\R_{>0}$-bundle over $M$. 

Let $G$ be the holonomy groupoid of $F$.  We let $r, s\colon G \rightarrow M$ denote
the range and source maps of $G$.  When $x$ is a point of $M$, the set of holonomies with source $x$ is denoted by $G_x$.  When $x$ and $y$ are points on the same leaf of $F$, we put $G_x^y = \ensemble{\gamma \in G \suchthat s(\gamma) = x, r(\gamma) = y}$.  The composition of holonomies are written as
\[
G_y^z \times G_x^y \rightarrow G_x^z, \quad (\gamma, \gamma') \mapsto \gamma \gamma'
\]
for any triple $x, y, z$ of points of $M$ on the same leaf.  When $T \subset M$ is a transversal (of dimension $q$) of $F$, the restricted
groupoid $G_T = \ensemble{ \gamma \in G \suchthat r\gamma, s\gamma \in T}$
becomes an \'{e}tale groupoid over $T$.

Even though $G_T$ might not be a Hausdorff
space, it is possible to take a covering of $G_T$ by open sets that are homeomorphic to $\R^q$.  A function on $G_T$ is said to be compactly supported smooth if it is so when restricted to one such open set.  The space of compactly supported smooth functions is denoted by $C_c^\infty(G_T)$.  As in~\cite{MR679730}*{Section~6}, the formula $f * f'(\gamma) = \sum_{\gamma = \gamma' \gamma''} f(\gamma') f'(\gamma'')$ defines an associative algebra structure on $C_c^\infty(G_T)$.

There is a naturally induced action of
the groupoid $G$ on the bundle $\tau$ over $M$.  We make an identification of the tangent bundle $\Tang T$ of $T$ with $\tau|_T$
given by the composition of the inclusion $\Tang T \rightarrow \Tang M |_T$ and
the projection $\Tang M |_T \rightarrow \tau|_T$.  The action of $G$ on $\tau$
induces a one of $G_T$ on $\Tang T$ via this identification.

A $G_T$-quasiinvariant measure on $T$ give rise to a positive definite
functional over the convolution algebra of $G_T$.  In this paper we concentrate
on the ones of the Lebesgue absolute continuity class among the transverse densities.  Hence
they are given by the sections of the $\R_{>0}$-bundle $\Det(T)$ which are Lebesgue
measurable and two such sections are identified when they agree off a negligible subset.

Given a foliation $F$ on $M$, a transversal $T$, and a density $\omega$ on $T$,
one obtains an inner product on $C^\infty_c(G_T)$ by
\[
(f_0, f_1) = \int_{x \in T} (f_0, f_1)_{\ell^2 G_x} \omega_x = \int_{x \in T} \sum_{\gamma \in G_x} f_0(\gamma) \overline{f_1(\gamma)} \omega_x.
\]
Let $H$ denote the Hilbert space completion of $C^\infty_c(G_T)$ with respect to
this inner product.  Let us recall the definition of relevant operator
algebras~\citelist{\cite{MR548112}\cite{MR679730}} on $H$ associated with the
groupoid $G_T$.

Let $\smooth{A}$ denote the $*$-algebra of the functions in $C^\infty_c(G_T)$
endowed with the convolution product and the involution $f^*(\gamma) =
\overline{f(\gamma^{-1})}$.  The left convolution
\[
\pi_l(f) \xi (\gamma) = \sum_{\gamma = \gamma' \gamma''} f(\gamma') \xi(\gamma'')
\]
defines a $*$-representation of $\smooth{A}$ on $H$.  The operator
norm closure $C^*_r(G_T)$ of $\smooth{A}$ is called the reduced C$^*$-algebra of
the foliation $(M; F)$.  The weak closure $W(M; F)$ of $\pi_l(\smooth{A})$ is called the von Neumann algebra of $(M; F)$.  These completions contain appropriate enlargements of $C^\infty_c(G_T)$ such as the space of compactly supported continuous functions $C_c(G_T)$, and we make use of these implicitly when there is no confusion.

For each leaf $l$ of $F$, let $\tilde{l}$ be the holonomy covering of $l$, and
$\pi$ be the covering map $\tilde{l} \rightarrow l$.  Thus, $\tilde{l}$ can be regarded as a copy of $G_x$ for any $x \in l$.  Then $W(M; F)$ can be described as the algebra of families of
operators $(T_l)_{l \in M/F}$ indexed by the space of the leaves of $F$ where $T_l \in B(\ell^2 (\pi^{-1}T \cap \tilde{l}))$, bounded in the sense that there exists a positive number $C$ satisfying $\norm{T_l} < C$ for all $l$, and measurable in the sense that the function $(\gamma, \gamma') \mapsto
(T_l e_{\gamma'}, e_{\gamma \gamma'})$ is measurable on $G_T \times_T G_T$, and two presentations
$(T_l)_{l \in M/F}$ and $(T_l')_{l \in M/F}$ are identified if they agree on a union of leaves
whose complement has measure $0$.

The commutant $W(\nu_\omega)$ of $\pi_l(\smooth{A})$ on $H$ is weakly generated by the
right convolution operators.  Namely, there is an anti-representation of $\smooth{A}$
on $H$ by
\[
\pi_r(f) \xi (\gamma) = \sum_{\gamma = \gamma' \gamma''} \xi(\gamma') f(\gamma'').
\]
Contrary to the case of left convolution, this representation may not be a
$*$-representation, and this failure can be measured by the modular function $\delta^{(\omega)}$ defined below.  Nevertheless, the operators of the form $\pi_r(f)$ for $f \in C_0(G_T)$ are a bounded operator in the
commutant of $\pi_l(\smooth{A})$.

The center $Z(W(M; F)) = W(M; F) \cap W(\nu_\omega)$ of $W(M; F)$ is naturally
identified with the algebra $L^\infty(T)^{G_T}$ consisting of bounded measurable functions on
$T$ which are invariant under the action of $G_T$.  Hence $W(M; F)$ is a factor
if and only if the action of $G_T$ on $T$ is ergodic with respect to the Lebesgue
measure class.

The density $\omega$ also gives a family $a_l$ of operator-valued
densities~\citelist{\cite{MR636521}\cite{MR1303779}}.  When $X$ is a section of
$\wedge^q \tau |_T = \wedge^q \Tang T$ and $l \in M/F$, the corresponding
operator $a(X)_l$ on $H_l = \ell^2 (\pi^{-1} T \cap \tilde{l})$ is given by the diagonal operator
$e_x \mapsto \pairing{\omega_x}{X_x} e_x$.  Then it defines a weight
$\phi^{(\omega)}$ on $W(M; F)$ by the invariant integral
\[
\phi^{(\omega)}((T_l)_{l \in M/F}) = \int_T (a_l T_l e_x, e_x) =
\int_T (T_l e_x, e_x) \omega_x
\]
on $T$.  The representation $\pi_l$ of $\smooth{A}$ on $H$ can be identified with the GNS representation of $\phi^{(\omega)}$.

The modular function $\delta^{(\omega)}$ from $G_T$ to $\R_{>0}$ associated with $\omega$ is defined by
\begin{equation}
  \label{eq:defn-modular-func}
  \omega_{s\gamma} = \delta^{(\omega)}(\gamma) \gamma^* \omega_{r \gamma}.
\end{equation}
Then $\omega$ is invariant under the holonomy transformations if and only if $\delta^{(\omega)}$ is identically equal to $1$ on $G_T$.

\begin{dfn}
The transverse density $\omega$ is said to be \textit{projectively invariant} if the function $\delta^{(\omega)}$ is locally constant on $G_T$. 
\end{dfn}

Suppose that $\omega$ is projectively invariant, and let $\gamma \in G_T$ be represented by the germ of a local
  holonomy map $g\colon S \rightarrow g(S)$ defined on an open connected
  neighborhood $S$ of $s(\gamma)$.  Then the set $U_g$ of the holonomies defined
  by $g$ becomes an open connected neighborhood of $\gamma$ in $G_T$ and $\delta^{(\omega)}$ becomes a constant function on $U_g$.  Thus, the pullback of $\omega|_{g(S)}$ by $g$ is a scalar multiple of $\omega|_S$ by the positive number $\delta^{(\omega)}(\gamma)$.

Let $\sigma^{(\omega)}_t$ be the modular automorphism associated with the weight
$\phi^{(\omega)}$.  Thus, when $f$ is a function of compact support on $G_T$, its effect on $\pi_l(f)$ is given by
\begin{equation}
\label{eq:right-convol-mod-auto}
\sigma^{(\omega)}_t(\pi_l(f)) = \pi_l((\delta^{(\omega)})^{i t} f),
\end{equation}
where $(\delta^{(\omega)})^{i t} f$ denotes the pointwise product of the
functions $\delta^{(\omega)}(\gamma)^{i t}$ and $f(\gamma)$ on $G_T$.

Suppose that $\omega$ is smooth.  Then $\delta^{(\omega)}$ is smooth on $G_T$, and~\eqref{eq:right-convol-mod-auto} implies that the $1$-parameter group $\sigma^{(\omega)}_t$ preserves the
subalgebra $\pi_l(C^{\infty}_c(G_T))$.  When this is the case, we let
$\R \ltimes_{\sigma^{(\omega)}} C^\infty_c(G_T)$ denote the linear space
$\schwfuncs(\R; C^\infty_c(G_T))$ of the Schwartz class functions on $\R$ with
values in $C^\infty_c(G_T)$, endowed with the convolution product twisted by the
action $\sigma^{(\omega)}$.

Let $\Det (T)$ denote the restriction of $\Det(M)$ to
$T$.  By means of the
section $\omega$ of $\Det(T)$, it can be identified with the direct product $\R_{>0} \times T$.  The natural action of $G_T$ on
$\Det(T)$ is identified with the one
\[
\gamma . (r, x) = (\delta^{(\omega)}(\gamma) r, x)
\]
of $G_T$ on $\R_{>0} \times T$.  We let $\Det(G_T)$ denote the corresponding
groupoid $G_T \ltimes \Det(T)$.

There is an injective homomorphism
\begin{equation}
  \label{eq:det-bundle-crossed-prod-emb}
  C^\infty_c(\Det(G_T)) \rightarrow \R \ltimes_{\sigma^{(\omega)}} C^\infty_c(G_T)
\end{equation}
with dense image, given by the Fourier transform on the $\R$-coordinate in $\R \times T \simeq \Det(G_T)^{(0)}$.

\subsection{Projective invariance of transverse density}
\label{sec:proj-invar-transv}

In this paper we consider the foliations satisfying the following two
conditions.  The first one is the following:
\begin{description}
\item[\condprojinv] There exists a nonvanishing smooth density $\omega \in
  \Gamma(T, \Det(T))$ which is projectively invariant under $G_T$.
\end{description}

When the groupoid $G_T$ satisfies the condition \condprojinv, the subgroupoid
$G_T^{(u)} = \ensemble{ \gamma \in G_T \suchthat \delta^{(\omega)}(\gamma) = 1}$
of $G_T$ is an open subset of $G_T$.  The second condition concerns the ergodicity of this subgroupoid.

\begin{description}
\item[\conduniterg] The action of $G_T^{(u)}$ is ergodic with respect to the
  Lebesgue measure class on $T$.
\end{description}

\begin{rmk}
  \label{rmk:ue-ergod}
  If $G_T$ and $\omega$ satisfy the condition {\condprojinv}, so do $G_S$ and
  $\omega|_S$ for any open set $S$ of $T$.  Similarly if they satisfy
  {\conduniterg}, the groupoid $G_S$ is Morita equivalent to $G_T$ by the
  groupoid bimodule
  \[
  G^S_T = \ensemble{ \gamma \in G \suchthat r(\gamma) \in S, s(\gamma) \in T},
  \]
  endowed with a left action of $G_S$ and a right action of $G_T$.  An analogous
  statement holds for $G_S^{(u)}$ and $G_T^{(u)}$.
\end{rmk}

\begin{rmk}
  By imposing the condition {\conduniterg}, we exclude the
  following kind of `false' examples.  When $M = \R$ and $F$ is given by the $\Tang M$ itself, the
  countable subset $T = \Z$ of $M$ is a transversal for $F$.  Then a transverse
  density over $T$ is equivalent to a sequence of numbers indexed by integers.
  Hence any choice of density on $T$ satisfy the assumption {\condprojinv}.  The
  associated foliation algebra is the convolution algebra in one real variable,
  which should be regarded as a trivial object in the category of operator algebras.
\end{rmk}

\begin{rmk}
  As a consequence of {\conduniterg}, we obtain that $C^*_r(G_T)$ and
  $C^*_r(G_T^{(u)})$ are simple and that the relative commutant of $C^*_r(G_T^{(u)})$ in $W(M; F)$ is trivial.
\end{rmk}

Given a transversal $T$, let $X_T$ denote the closure of the image of the image
of the module morphism $\delta^{(\omega)}\colon G_T\rightarrow \R_{>0}$.  We
have
\[
X_T = \Sp \left (\sigma^{(\omega)} \right )
\]
by~\eqref{eq:right-convol-mod-auto}.  Let $S(W(M; F))$ denote the $S$-invariant
of $W(M; F)$~\cite{MR0303306}.

\begin{lem}
  \label{lem:proj-inv-density-s-set-calc}
  Suppose that $G_T$ and $\omega$ satisfy the conditions {\condprojinv} and
  {\conduniterg}.  Then $S(W(M; F)) \cap \R_{>0}$ is equal to $X_T$.
\end{lem}

\begin{proof}
  The centralizer algebra $W(M; F)^{\phi_\omega}$ is identified with the groupoid
  von Neumann algebra $W(G_T^{(u)})$, which is a factor by the condition
  {\conduniterg}.  Then one has $\Gamma(\sigma^{(\omega)}) =
  \Sp \left (\sigma^{(\omega)} \right ) = X_T$ by~\cite{MR0303306}*{Proposition~2.2.2 (c)},
  which proves the assertion.
\end{proof}

\begin{example}
  A foliation $F$ is said to be transversely affine when there exists a covering
  of $M$ by some foliation charts with respect to which the holonomy maps become affine transformations.  Suppose that $T$ is a
  transversal contained in a foliation chart
  \[
  U \simeq \ensemble{ (x_1, \ldots, x_q, y_1, \ldots, y_p) \in \R^p \times \R^q}
  \]
  of $F$ such that the holonomy transformations are affine with respect to the
  transverse coordinate system $(x_1, \ldots, x_q)$.  Then the transverse volume
  form
  \[
  \absolute{d x} = \absolute{d x_1 \wedge \cdots \wedge d x_q}
  \]
  gives a projectively invariant density on $T$.
\end{example}

\begin{example}
\label{exmpl:hyperb-torus}
  Let $A$ be a hyperbolic matrix in $\SL_2(\Z)$.  It defines an homeomorphism of
  $\T^2 = \R^2/\Z^2$.  The associated mapping cone $M = M_A \T^2$ admits an
  Anosov foliation $F$.  Specifically, let $u$ and $v$ be eigenvectors of $A$,
  associated with eigenvalues $\lambda < 1 < \lambda^{-1}$.  Then the subspace of
  $\Tang M$ spanned by the image of $u$ and the suspension flow is integrable.
  We obtain the hyperfinite factor of type III$_\lambda$ as the associated von
  Neumann algebra~\cite{MR0461584}.
\end{example}

\begin{example}
  Let $X$ be the cosphere bundle $S^* \Sigma_g$ of an oriented closed Riemannian
  surface of genus $g > 1$.  Then $X$ admits a foliation $F$ of dimension $2$
  generated by the holocyclic flow and the geodesic flow.  If one takes an
  embedding of $\pi_1(\Sigma_g)$ into $SL_2(\R)$ and identify $X$ with $SL_2(\R)
  / \pi_1(\Sigma_g)$, $F$ is generated by the actions the following matrices
  from left:
  \begin{align*}
    A_t &= \left[ \begin{array}{cc} 1 & t \\ 0 & 1 \end{array} \right], &
    B_t &= \left[ \begin{array}{cc} e^{t} & 0 \\ 0 & e^{-t} \end{array} \right].
  \end{align*}
  When $T$ is transversal to $F$, the multiplication from left by
  \begin{equation*}
    C_t = \left[ \begin{array}{cc} 1 & 0 \\ t & 1 \end{array} \right]
  \end{equation*}
  gives a projectively invariant coordinate on $T$.
\end{example}

\subsection{Transverse Metric Trivialization}
\label{sec:transv-metr-triv}

By assuming the condition {\condprojinv}, we are limiting our consideration to a particular class of foliations with some tame transverse measure structure, but such a foliation does not need to possess a projectively invariant
transverse metric structure.  The lack of such a metric is an obstacle to the
problem of extending cyclic cocycles on the smooth convolution algebra to the
C$^*$-algebra of $G_{T}$.  To remedy this, we introduce the groupoid
$\Met(G_{T})$ of transverse metric trivialization and work with the cyclic
cocycles defined on its groupoid algebra as in~\cite{MR866491}.

For each $x \in M$, let $\Met_{x}$ denote the set of the strictly positive
quadratic forms on the vector space $\tau_{x}$.  By functoriality there is a natural action of $G_T$ on the manifold $\Met(T) = \bigcup_{x \in T}
\Met_{x}$.  Let $\Met(G_{T})$ denote the associated
groupoid.

The module map~\eqref{eq:defn-modular-func} induces a
groupoid homomorphism $\delta^{(\omega)}$ from $\Met(G_{T})$ to $\R_{>0}$
by the composition of the natural projection $\Met(G_{T}) \rightarrow
G_{T}$ and $\delta^{(\omega)}$ on $G_T$.

On one hand, any density $\omega$ defines a subset of the total space of
$\Det(M)$ as its image.  On the other hand, there is a natural surjection
$\Met(M) \rightarrow \Det(M)$.  We let $\SMet(M)$ denote the inverse image of
$\omega$ inside $\Met(M)$.  When $\gamma \in G_T$ and $\xi \in \SMet(T)_{s
  \gamma}$, the element
\[
\SMet(\gamma)\xi = \delta^{(\omega)}(\gamma)^{-1} \Met(\gamma) \xi
\]
is in $\SMet(T)_{r \gamma}$.  It follows that $\gamma \mapsto \SMet(\gamma)$
defines an action of $G_T$ on $\SMet(T)$.  Let us denote $\SMet(G_T) = G_T
\ltimes_{\SMet} \SMet(T)$.  Since the fiber of the projection $\SMet(T)
\rightarrow T$ admits a canonical spin structure, one has the natural
isomorphism between the $K$-groups of $C^*_r(\SMet(G_T))$ and $C^*_r(G_T)$.

The modular function $\delta^{(\omega)}$ on $G_T$ lifts to a function (again
denoted by $\delta^{(\omega)}$ by abuse of notation) on
$\SMet(G_T)$.  It defines a one-parameter group
$\sigma^{(\omega)}$ on $C^\infty_c(\SMet G_T)$.  We have an
embedding
\begin{equation}
  \label{eq:met-bundle-emb-cross-smet-bundle}
  C^\infty_c(\Met G_T) \rightarrow
  \R \ltimes_{\sigma^{(\omega)}} C^\infty_c (\SMet G_T)
\end{equation}
of dense image into the smooth crossed product analogous to
\eqref{eq:det-bundle-crossed-prod-emb}.

\section{Periodic action and cyclic cocycles}
\label{sec:peri-acti-cycl}

Suppose that $\sigma$ is an action of $\T$ on a C$^*$-algebra $B$.  For each $n \in \Z$, let $B_n$ denote the spectral subspace $\ensemble{ x \in B \suchthat \sigma_t(x) = e^{2 \pi i n t} x}$.  We let $A$ denote the fixed point algebra $B_0$ of $\sigma$.

The dual action $\hat{\sigma}$ of $\Z$ on $\T \ltimes_\sigma B$ induces the
Takesaki--Takai duality
\[
\Z \ltimes_{\hat{\sigma}} \T \ltimes_\sigma B \simeq_{\KK} B.
\]
The Pimsner--Voiculescu exact sequence~\cite{MR587369} for the automorphism $\hat{\sigma}$ of $\T \ltimes_\sigma B$ is given by the $6$-term exact
sequence
\begin{equation}
  \label{eq:pimsner-voiculescu-6-term-seq}
  \begin{CD}
    K_0(\T \ltimes_\sigma B) @>{\Id - \hat{\sigma}_*}>> K_0(\T \ltimes_\sigma B) @>{\iota_*}>> K_0(B) \\
    @A{\partial}AA & & @VV{\partial}V \\
    K_1(B) @<<{\iota_*}< K_1(\T \ltimes_\sigma B) @<<{\Id - \hat{\sigma}_*}< K_1(\T \ltimes_\sigma B)
  \end{CD},
\end{equation}
where $\iota^*$ is induced by the inclusion homomorphism $\iota\colon \T \ltimes_\sigma B \rightarrow \Z \ltimes_{\hat{\sigma}} \T \ltimes_\sigma B$ and the natural isomorphism $K_*(\Z \ltimes_{\hat{\sigma}} \T \ltimes_\sigma B) \simeq K_*(B)$.

The action $\sigma$ may be regarded as an action of $\R$ via the surjection
$\R \rightarrow \R / \Z \simeq \T$.  Let $\Psi$ be the $*$-algebra homomorphism from
$\R \ltimes_\sigma B$ to $\T \ltimes_\sigma B$ characterized as the unique
extension of the mapping
\begin{equation}
  \label{eq:mapping-torus-periodization-defn}
  \Psi(f)_t = \sum_{m \in \Z} f(t + m)
\end{equation}
from $L^1(\R; B) \subset \R \ltimes_\sigma B$ to $L^1(\T; B) \subset \T
\ltimes_\sigma B$.

\begin{lem}
  \label{lem:pv-boundary-res-to-dual-of-torus}
  The map $\Psi_*\colon K_* (\R \ltimes_\sigma B) \rightarrow K_*(\T \ltimes_\sigma B)$
  induced by $\Psi$ is identified with the map $\partial$ in~\eqref{eq:discr-spec-R-action-fixed-point-PV-seq}, via the Connes--Thom
  isomorphism $K_* (B) \simeq K_{*+1}( \R \ltimes_\sigma B)$ and the natural
  isomorphism $K_*(\T \ltimes_\sigma B) \simeq K_*(A)$ induced by~\eqref{eq:mor-equiv-fixed-t-action-from-discr-spec-r-action}.
\end{lem}

\begin{proof}
  The crossed product $\R \ltimes_\sigma B$ is identified with the mapping torus
  \[
  M_{\hat{\sigma}}(\T \ltimes_\sigma B) = \ensemble{ f \in C(\R; \T \ltimes_\sigma B) \suchthat f_{t + 1} = \hat{\sigma}(f_t) }.
  \]
  Under this identification the homomorphism $\Psi$ corresponds to the
  evaluation map $f \mapsto f_0$ from $M_{\hat{\sigma}}(\T \ltimes_\sigma B)$ to
  $\T \ltimes_\sigma B$.  Then Connes's proof~\cite{MR605351}*{pp. 48--49} of the
  Pimsner--Voiculescu exact sequence implies the assertion.
\end{proof}

Let $\smooth{B}$ be a locally convex algebra over $\C$.  Let $\Omega(\smooth{B})$ be the
universal differential graded algebra over $\smooth{B}$.  It
is the direct sum $\oplus_{n \in \N} \Omega_n(\smooth{B})$, where
\[
\Omega_0(\smooth{B}) = \smooth{B}, \quad \Omega_n(\smooth{B}) = \smooth{B}^{\otimes
  n} \oplus \smooth{B}^{\otimes n + 1}\quad (n > 0).
\]
An element $(a^1, \ldots, a^n) \oplus (b^0, \ldots, b^n) \in \Omega_n(\smooth{B})$ is understood to represent the $n$-form $d a^1 \cdots d a^n + b^0 d b^1 \cdots d b^n$.  Thus the  differential $d\colon \Omega_n(\smooth{B}) \rightarrow
\Omega_{n+1}(\smooth{B})$ is given by the identity map on the factor
$\smooth{B}^{\otimes n + 1}$ in both sides.  The algebra structure on $\Omega(\smooth{B})$ can be determined using the Leibniz rule for $d$.

Let $\kappa$ denote the operator
\[
f^0 \otimes \cdots \otimes f^{n} \mapsto (-1)^n f^n \otimes f^0 \otimes f^1 \otimes \cdots \otimes f^{n-1}
\]
on $\smooth{B}^{\otimes n+1}$.  It satisfies $\kappa^{n+1} = 1$.  Next, there is
the $b'$-operator
\[
b' \psi(f^0, \ldots, f^{n}) = \sum_{k=0}^{n-1} (-1)^{k} \psi(f^0, \ldots, f^k f^{k+1}, \ldots, f^{n})
\]
and the Hochschild coboundary operator $b$
\[
b\psi(f^0, \ldots, f^{n}) = b' \psi(f^0, \ldots, f^{n}) + (-1)^{n} \psi(f^{n}f^0, f^1, \ldots, f^{n-1}),
\]
from $(\smooth{B}^{\otimes n})'$ to $(\smooth{B}^{\otimes n+1})'$.  The cyclic
cocycles are precisely the ones in the joint kernel of $1 - \kappa$ and $b$.
One also has the $S$-operator of Connes~\cite{MR823176}*{Section I.4}.

As in the C$^*$-algebraic setting before, let $\sigma$ be an action of $\T$ on $\smooth{B}$.  In the rest of this
section we let $\phi$ denote a $\sigma$-invariant cyclic $n$-cocycle on
$\smooth{B}$.  One obtains the
following two new cyclic cocycles.

We let $\R \ltimes_\sigma \smooth{B}$ denote the space of Schwartz class functions from $\R$ to $\smooth{B}$ endowed with the convolution product.  The first cocycle, $\hat{\phi}$ is a cyclic $n$-cocycle on $\R \ltimes_\sigma \smooth{B}$ given
by
\[
\pairing{\hat{\phi}}{f^0 d f^1 \cdots d f^n} = \int_{\sum_{j=0}^n t_j = 0} \phi(f^0_{t_0}, \sigma_{t_0}(f^1_{t_1}), \cdots,  \sigma_{\sum_{j=0}^{n-1}t_j}(f^n_{t_n}).
\]
Let $D$ be the generator of the action $\sigma$.  The second, $i_D \phi$ is a cyclic $n + 1$-cocycle on $\smooth{B}$ defined by
\[
i_D \phi(f^0, \cdots, f^{q' + 1}) = \sum_{j=1}^{q'} (-1)^{j-1}
\pairing{\phi}{f^0 d f^1 \cdots D(f^j) \cdots d f^{q' + 1}}.
\]

There is also a boundedness condition for $\phi$, namely the $n$-trace property, which depends on the choice of a norm $\bnorm{\cdot}$ on $\smooth{B}$.  Specifically, $\phi$ is said to be an $n$-trace if for any elements $(g^i)_{i = 1}^n$ of $B$, there exists a constant $C$ satisfying
\[
\absolute{\tau_\phi(f^0 g^1d f^1 \cdots g^n d f^n)} \le C \prod_{i = 1}^n \bnorm{g^i} \quad (f^i \in \smooth{B})
\]
where $\tau_\phi$ is the associated closed graded trace on $\Omega^*(\smooth{B})$.  If $\phi$ satisfies the $n$-trace condition, the map $K_n(\smooth{B}) \rightarrow \C$ induced by $\phi$ can be extended to $K_n(\banach{B})$, where $\banach{B}$ is the completion of $\smooth{B}$ with respect to $\bnorm{\cdot}$.  Moreover, when this is the case, $\hat{\phi}$ and $i_D \phi$ are also $n$-trace and $(n+1)$-trace on $\R \ltimes
\smooth{B}$ and $\smooth{B}$ respectively~\cite{MR1303779}*{Section
  3.6.$\beta$}.  When the norm $\bnorm{ \cdot }$ is chosen appropriately, the algebra $\banach{B}$ becomes a spectral subalgebra of the initial C$^*$-algebra $B$ and have the same $K$-groups.
  
The action $\sigma$ on $\smooth{B}$ has a canonical extension to the algebra
$\Omega(\smooth{B})$.  Then the closed trace $\tau_\phi$ on $\Omega(\smooth{B})$
associated with $\phi$ becomes invariant under this action.  Given an element $f$
of $\R \ltimes_\sigma \smooth{B}$, we let $d f$ denote the ``pointwise
derivation'' $t \mapsto d (f_t)$, regarded as an element in $\R \ltimes_\sigma
\Omega(\smooth{B})$.

For the ease of notation, we introduce several notations which are only used in
this section.  Given the standard coordinate $(t_0, \ldots, t_n)$ on $\R^{n+1}$
and an integer $1 \le j < n$, let $s_j$ denote the quantity $\sum_{0 \le k < j}
t_k$.  In addition, given the elements $f^0, \ldots, f^n$ of $\R \ltimes_\sigma
\smooth{B}$, let $g_j$ denote the function $\sigma_{s_j}(f^j_{t_j})$ from
$\R^{j}$ to $\smooth{B}$ for $0 \le j \le n$.

\begin{lem}
  For each $m \in \Z$, the functional
  \[
  \hat{\phi}_m(f^0,\ldots, f^n) = \tau_\phi((f^0 d f^1 \cdots d f^n)_m) = \int_{(t_j) \in \R^n, \sum_{j=0}^n t_j = m} \phi(g_0, \ldots, g_n)
  \]
  is a Hochschild $n$-cocycle over $\R \ltimes_\sigma \smooth{B}$.
\end{lem}

\begin{proof}
  By the $\sigma$-invariance of $\phi$ and the assumption that $\sigma$ is
  periodic implies
  \[
  \phi(f^{n+1}_{t_{n+1}} \sigma_{t_{n+1}}(g_0), \ldots, \sigma_{t_{n+1}}(g_n)) = \phi(g_{n+1} g_0, g_1, \ldots, g_n)
  \]
  whenever one has $\sum_{j=0}^{n+1} t_j = m$.  Thus, $b \hat{\phi}_m(f^0,
  \ldots, f^{n+1})$ is equal to
  \[
  \int_{\sum_{j=0}^{n+1} t_j = m} \sum_{k=0}^{n} (-1)^{k} \phi(g_0, g_1, \ldots, g_j g_{j+1}, \ldots, g_n) + (-1)^{n+1} \phi(g_{n+1} g_0, g_1, \ldots, g_n).
  \]
  By the Hochschild cocycle condition on $\phi$, the integrand of the above
  formula vanishes.  Hence one has $b \hat{\phi}_m = 0$.
\end{proof}

\begin{rmk}
  The Hochschild cocycle $\hat{\phi}_m$ satisfies the cyclicity condition when
  $m = 0$, but otherwise there is no such guarantee.
\end{rmk}

Let $\hat{D}$ be the derivation $\hat{D}(f)_t = t f_t$ on $\R \ltimes
\smooth{B}$.  It is the generator of the dual action $\hat{\sigma}$ of $\R$ on
$\R \ltimes_\sigma \smooth{B}$.

\begin{lem}
  \label{lem:eta-m-is-b-cocycle}
  For each $m \in \Z$, the multilinear functional
  \[
  \eta_m(f^0, \ldots, f^{n+1}) = \sum_{k=1}^{n+1}(-1)^{k+1} \tau_\phi((f^0 d f^1 \cdots \hat{D}(f^k) \cdots d f^{n+1})_m)
  \]
  is an $(n+1)$-Hochschild cocycle over $\R \ltimes_\sigma \smooth{B}$.
\end{lem}

\begin{proof}
  By definition, $\eta_m$ is the action $\hat{D} \# \hat{\phi}_m$ of the
  derivation $\hat{D}$ on the Hochschild cocycle
  $\hat{\phi}_m$~\cite{MR1303779}*{Remark~3.2.30.b}.  Hence it satisfies the
  Hochschild coboundary relation.
\end{proof}

\begin{lem}
  \label{lem:b-prime-psi-m-in-image-of-id-minus-karoubi-op}
  For any integers $n$ and $m$, one has $(1-\kappa^*)\eta_m = (-1)^n m b' \hat{\phi}_m$.
\end{lem}

\begin{proof}
  For each $1 \le k \le n$, the term $\tau_\phi((f^0 d f^1 \cdots \hat{D}(f^k)
  \cdots d f^{n+1})_m)$ contributes as
  \[
  \int_{t_0 + \cdots + t_n = m} t_k \sum_{j=0}^{k-1} (-1)^{k - j - 1} \phi(g_0, \ldots, g_j g_{j+1}, \ldots, g_n).
  \]
  Hence $\eta_m(f^0, \ldots, f^n)$ is the sum of these terms.  Similarly,
  $\kappa^*\eta_m(f^0, \ldots, f^n)$ is the sum of
  \begin{multline*}
    (-1)^n \int_{t_0 + \cdots + t_n = m} t_k \phi(f^n_{t_n} \sigma_{t_n}(g_0), \sigma_{t_n}(g_1), \ldots, \sigma_{t_n}(g_{n-1}))\\
    + t_k \sum_{j=0}^{k-1} (-1)^{k - j - 1} \phi(f^n_{t_n}, \sigma_{t_n}(g_0), \ldots, \sigma_{t_n}(g_j) \sigma_{t_n}(g_{j+1}), \ldots, \sigma_{t_n}(g_{n-1}))
  \end{multline*}
  for $0 \le k \le n - 1$.

  First, by the cyclicity and the $\sigma$-invariance of $\phi$, one has
 \begin{multline*}
 (-1)^n (-1)^{k - j - 1} \phi(f^n_{t_n}, \sigma_{t_n}(g_0), \ldots, \sigma_{t_n}(g_j) \sigma_{t_n }(g_{j+1}), \ldots, \sigma_{t_n}(g_{n-1})) \\
 = (-1)^{k - j - 1} \phi(g_0, \ldots, g_j g_{j+1}, \ldots, g_n)
 \end{multline*}
 for each $j$ when $t_0 + \cdots + t_n = m$.  This shows that
 \begin{multline*}
 \sum_{j=0}^{k-1} (-1)^{k - j - 1} \phi(g_0, \ldots, g_j g_{j+1}, \ldots, g_n)\\
 - (-1)^n \sum_{j=0}^{k-1} (-1)^{k - j - 1} \phi(f^n_{t_n}, \sigma_{t_n}(g_0), \ldots, \sigma_{t_n}(g_j) \sigma_{t_n}(g_{j+1}), \ldots, \sigma_{t_n}(g_{n-1}))
 \end{multline*}
 vanishes for $1 \le k \le n - 1$.

 Next, by the Hochschild cocycle condition on $\phi$, one has
\begin{multline*}
  t_n \sum_{j=0}^{n-1} (-1)^{n - j - 1} \phi(g_0, \ldots, g_j g_{j+1}, \ldots, g_n)\\
  = (-1)^{n+1} t_n \phi(f^n_{t_n}, \sigma_{t_n}(g_0), \ldots, \sigma_{t_n}(g_j) \sigma_{t_n}(g_{j+1}), \ldots, \sigma_{t_n}(g_{n-1})).
  \end{multline*}
  Combining these, $(1 - \kappa^*)\eta_m(f^0, \ldots, f^n)$ reduces to
 \[
 \int_{\sum_{j=0}^3 t_j = m} (\sum_{j=0}^n t_j) \phi(f^n_{t_n} \sigma_{t_n}(f^0_{t_0}), \sigma_{t_n + s_1}(f^1_{t_1}), \ldots, \sigma_{t_n + s_{n-1}}(f^{n-1}_{t_{n-1}})),
 \]
 which proves the assertion in this case.
\end{proof}

Let us illustrate the the content of Lemma~\ref{lem:b-prime-psi-m-in-image-of-id-minus-karoubi-op} for the case $n = 1$.  By definition of $\eta_m$, one has
\begin{equation*}
\begin{split}
\eta_m(f^0, f^1, f^2) &= \pairing{\tau_\phi}{\left(f^0 \hat{D}(f^1) d f^2 - f^0 d f^1 \hat{D}(f^2)\right)_m}\\
&= \pairing{\tau_\phi}{\left(f^0 \hat{D}(f^1) d f^2 - f^0 d (f^1 \hat{D}(f^2)) + f^0 f^1 d \hat{D}(f^2)\right)_m},
\end{split}
\end{equation*}
which is equal to
\begin{multline*}
\int_{t_0 + t_1 + t_2 = m} t_1 \phi(f^0_{t_0} \sigma_{t_0}(f^1_{t_1}), \sigma_{t_0 + t_1}(f^2_{t_2})) \\
- t_2 \phi(f^0_{t_0}, \sigma_{t_0}(f^1_{t_1}) \sigma_{t_0 + t_1}(f^2_{t_2})) + t_2 \phi(f^0_{t_0} \sigma_{t_0}(f^1_{t_1}), \sigma_{t_0 + t_1}(f^2_{t_2})).
\end{multline*}

Then
\[
(1 - \kappa^*) \eta_m(f^0, f^1, f^2) = \eta_m(f^0, f^1, f^2) - \eta_m(f^2, f^0, f^1)
\]
can be computed as
\begin{multline*}
\int_{t_0 + t_1 + t_2 = m}  t_1 \phi(f^0_{t_0} \sigma_{t_0}(f^1_{t_1}), \sigma_{t_0 + t_1}(f^2_{t_2})) - t_2 \phi(f^0_{t_0}, \sigma_{t_0}(f^1_{t_1}) \sigma_{t_0 + t_1}(f^2_{t_2})) \\
+ t_2 \phi(f^0_{t_0} \sigma_{t_0}(f^1_{t_1}), \sigma_{t_0 + t_1}(f^2_{t_2}))  - t_0 \phi(f^2_{t_2} \sigma_{t_2}(f^0_{t_0}), \sigma_{t_2 + t_0}(f^1_{t_1})) \\
+ t_1 \phi(f^2_{t_2}, \sigma_{t_2}(f^0_{t_0}) \sigma_{t_2 + t_0}(f^1_{t_1})) - t_1 \phi(f^2_{t_2} \sigma_{t_2}(f^0_{t_0}), \sigma_{t_2 + t_0}(f^1_{t_1})).
\end{multline*}
By the cyclicity condition $\phi = \kappa^* \phi$ and the $\sigma$-invariance on $\phi$, one has
\[
t_1 \phi(f^0_{t_0} \sigma_{t_0}(f^1_{t_1}), \sigma_{t_0 + t_1}(f^2_{t_2})) = - t_1 \phi(f^2_{t_2}, \sigma_{t_2}(f^0_{t_0}) \sigma_{t_2 + t_0}(f^1_{t_1})).
\]
On the other hand, the Hochschild cocycle condition $b\phi = 0$ implies
\begin{multline*}
t_2 \left\{ \phi(f^0_{t_0}, \sigma_{t_0}(f^1_{t_1}) \sigma_{t_0 + t_1}(f^2_{t_2})) - \phi(f^0_{t_0} \sigma_{t_0}(f^1_{t_1}), \sigma_{t_0 + t_1}(f^2_{t_2}))\right\}\\
= t_2 \phi(f^2_{t_2} \sigma_{t_2}(f^0_{t_0}), \sigma_{t_2 + t_0}(f^1_{t_1})).
\end{multline*}
Combining these, one obtains
\begin{equation*}
(1 - \kappa^*) \eta_m(f^0, f^1, f^2) = \int_{t_0 + t_1 + t_2 = m} (-m) \phi(f^2_{t_2} \sigma_{t_2}(f^0_{t_0}), \sigma_{t_2 + t_0}(f^1_{t_1})),
\end{equation*}
which implies $(1 - \kappa^*)\eta_m = -m b' \hat{\phi}_m$ when $n=1$.

Analogously to the $\R$-crossed product case, we have the dual cyclic cocycle $\hat{\phi}^\T$ on $\T \ltimes_\sigma \smooth{B}$ defined by
\[
  \hat{\phi}^\T(f^0,\ldots, f^n) = \tau_\phi((f^0 d f^1 \cdots d f^n)_m) = \int_{(t_j) \in \T^n, \sum_{j=0}^n t_j = m} \phi(g_0, \ldots, g_n).
\]
When $\smooth{B}$ is Fr\'{e}chet, the formula~\eqref{eq:mapping-torus-periodization-defn} defines a homomorphism from $\R \ltimes_{\sigma} \smooth{B}$ to $\T \ltimes_{\sigma} \smooth{B}$, which we denote again by $\Psi$.  Then the pullback cocycle on $\R \ltimes_\sigma \smooth{B}$ is represented as
 \begin{equation}
 \label{eq:pb-torus-dual-by-periodicization}
 \Psi^* (\hat{\phi}^\T)(f^0, \ldots, f^n) = \sum_{m \in \Z} \pairing{\tau_\phi}{(f^0 d f^1 \cdots d f^n)_m} = \sum_{m \in \Z} \hat{\phi}_m(f^0, \ldots, f^n).
 \end{equation}
  We note that this pullback do not agree with $\hat{\phi}_0$.  However, they define the same periodic cyclic cohomology class, as seen in the next theorem.

\begin{thm}
\label{thm:dual-cocycle-cohom-R-and-T}
Let $\sigma$ be an action of $\T$ on a Fr\'{e}chet algebra $\smooth{B}$.  When $\phi$ is a $\sigma$-invariant cyclic $n$-cocycle on $\smooth{B}$, the cyclic cocycles $\hat{\phi}$ and $\Psi^*(\hat{\phi}^\T)$ are cohomologous in $\HC^{n+2}(\R \ltimes_{\sigma} \smooth{B})$.
\end{thm}

\begin{proof}
  We know that $m S \hat{\phi}_m$ is a cyclic coboundary for any $m$ by Lemmas~\ref{lem:eta-m-is-b-cocycle} and~\ref{lem:b-prime-psi-m-in-image-of-id-minus-karoubi-op}, and the $(b, B)$-bicomplex
 argument as in~\citelist{\cite{MR823176}*{Lemma~29}\cite{MR2052770-Cuntz}}.  By~\eqref{eq:pb-torus-dual-by-periodicization}, $\Psi^*(\hat{\phi}^\T)$ becomes cohomologous to $\hat{\phi}_0$ after applying $S$.  Since one has
 $\hat{\phi}_0 = \hat{\phi}$, this proves the assertion.
\end{proof}

We now combine the considerations of the C$^*$-algebraic setting and the Fr\'{e}chet algebra setting to obtain the main result of this section.  Let $\sigma$ be an action of $\T$ on a C$^*$-algebra $B$.  Moreover, we assume that it acts
smoothly on a spectral subalgebra $\smooth{B}$ of $B$ with a Fr\'{e}chet topology.

\begin{cor}
\label{cor:period-dual-action-pairing-der-cycle-k-group}
  Let $\sigma$, $B$, and $\smooth{B}$ as above, and let $\phi$ be a $\sigma$-invariant $n$-trace on $\smooth{B}$.  Then, for any $x \in K_*(\smooth{B})$ one has
  $\pairing{i_D \phi}{x} = \pairing{\hat{\phi}^\T}{\partial(x)}$ with respect to the
  homomorphism $\partial$ in~\eqref{eq:pimsner-voiculescu-6-term-seq}.
\end{cor}

\begin{proof}
The cocycles $\hat{\phi}$ and $i_D \phi$ induce the same
  map~(implicitly proved in \cite{MR945014}*{Lemma~6.3}; see \cite{MR2738561}*{Proposition~14} for an explicit proof) via the identification given by
  Connes--Thom isomorphism $\Phi\colon K_*(\R \ltimes_\sigma B) \rightarrow
  K_{*+1}(B)$.  The assertion follows from Lemma~\ref{lem:pv-boundary-res-to-dual-of-torus} and Theorem~\ref{thm:dual-cocycle-cohom-R-and-T}.
\end{proof}

When $\sigma$ has the full strong spectrum, that is, if we have $B_{-n} B_n = A$ for any $n$, we have the strong Morita equivalence
\begin{equation}
  \label{eq:mor-equiv-fixed-t-action-from-discr-spec-r-action}
  \T \ltimes_\sigma B \simeq_{\KK} A.
\end{equation}
This condition is satisfied when $A$ is simple and $B_n \neq \ensemble{0}$ for every $n$, which is relevant to our setting in view of the condition {\conduniterg}.  When this is the case, the $6$-term exact sequence
\eqref{eq:pimsner-voiculescu-6-term-seq} for the action $\hat{\sigma}$ of $\Z$
on $\T \ltimes_\sigma B$ becomes
\begin{equation}
  \label{eq:discr-spec-R-action-fixed-point-PV-seq}
  \begin{CD}
    K_0(A) @>{\Id - \alpha_*}>> K_0(A) @>{\iota_*}>> K_0(B) \\
    @A{\partial}AA & & @VV{\partial}V \\
    K_1(B) @<<{\iota_*}< K_1(A) @<<{\Id - \alpha_*}< K_1(A).
  \end{CD}
\end{equation}

\begin{cor}
\label{cor:saturated-action-pimsner-voiculescu-nest}
 Let $\sigma$ be an action of $\R$ on a Fr\'{e}chet pre-C$^*$-algebra
  $\smooth{B}$.  Suppose that the spectrum of $\sigma$ is equal to $T \Z \subset
  \R$ for some real number $T$.  Then, for any $x \in K_*(\smooth{B})$ one has
  $\pairing{i_D \phi}{x} = T \pairing{\phi}{\partial(x)}$ with respect to the
  homomorphism $\partial$ in~\eqref{eq:mor-equiv-fixed-t-action-from-discr-spec-r-action}.
\end{cor}

\begin{proof}
   When $R$ is an arbitrary positive real number, we have $i_{R D} \phi = R i_D
  \phi$.  Hence we may assume that $T = 2 \pi$.  Thus we assume that $\sigma$
  comes from a full spectrum action of $\T = \R / 2 \pi \Z$ on $\smooth{B}$.
  
    The cyclic cocycle $\phi|_\smooth{A}$ corresponds to the multilinear functional
 \begin{equation*}
 \begin{split}
 \psi_0(f^0, \ldots, f^n) &= \pairing{\tau_\phi}{E((f^0 d f^1 \cdots d f^n)_0)}\\
 &= \int_T d t' \int_{t_0 + \cdots + t_n = 0} d t_0 \cdots d t_n \sigma_{t'}^*\phi\big(f^0_{t_0}, \sigma_{s_1}(f^1_{t_1}), \ldots, \sigma_{s_n}(f^n_{t_n})\big)\\
 &= \int_{t_0 + \cdots + t_n = 0} \phi\big(f^0_{t_0}, \sigma_{s_1}(f^1_{t_1}), \ldots, \sigma_{s_n}(f^n_{t_n})\big) d t_0 \cdots d t_n\\
  \end{split}
  \end{equation*}
  on $\schwfuncs(\T; \smooth{B}) \subset \T \ltimes_\sigma \smooth{B}$ via the
  strong Morita equivalence~\eqref{eq:mor-equiv-fixed-t-action-from-discr-spec-r-action}.
\end{proof}

\begin{example}
  \label{example:pv-nest-corr-torus-transl}
  Let $\smooth{B}$ be the Fr\'{e}chet algebra $C^\infty(\T)$ and $\sigma$ be the
  action of $\T$ on $\smooth{B}$ induced by the translation.  There is a
  unique $\sigma$-invariant trace $\tau$ on $\smooth{B}$, namely the integration
  of functions with respect to the normalized Haar measure on $\T$.  The fixed
  point algebra $\smooth{A}$ is equal to $\C$, and the cyclic $1$-cocycle $i_D
  \tau$ on $\smooth{B}$ is equal to the bilinear functional
 \[
 i_D \tau(f^0, f^1) = \frac{1}{2 \pi i} \int_\T f^0 d f^1,
 \]
 defined as the invariant integral of $1$-forms over the $1$-dimensional manifold $\T$.

 The $6$-term exact sequence~\eqref{eq:discr-spec-R-action-fixed-point-PV-seq}
 is the one induced by the extension of algebras $C_0(0, 1) \rightarrow C(\T)
 \rightarrow \C$.  The map induced by $\tau|_\C$ on $K_0(\C)$ and the one
 induced by $i_D \tau$ on $K_1(C_0(0, 1))$ indeed agree, as each of them sends
 generators of the corresponding $K$-group to $\pm 1$.
\end{example}

\begin{rmk}
  Let $\alpha$ be an action of $\Z$ on a C$^*$-algebra $A$. Suppose that there
  is a Fr\'{e}chet subalgebra $\smooth{A} \subset A$ which is invariant under
  $\alpha$.  R. Nest~\cite{MR961899} constructed a homomorphism $\partial^*\colon
  \HP^*(\smooth{A}) \rightarrow\HP^{*+1}(\Z \ltimes_\sigma \smooth{A})$ by means
  of a spectral sequence, which is transpose to the homomorphism $\partial$ in~\eqref{eq:pimsner-voiculescu-6-term-seq}.

  If $\phi$ is a $\alpha$-invariant cyclic $n$-cocycle on $\smooth{A}$, then one
  has the dual cocycle $\hat{\phi}^\Z$ on $\Z \ltimes_\alpha \smooth{A}$ by
  \[
  \hat{\phi}^\Z(f^0,\ldots, f^n) = \sum_{m_0 +\cdots + m_n = 0} \phi(f^0_{m_0}, \alpha_{m_0}( f^1_{m_1}), \ldots, \alpha_{m_0 + \cdots + m_{n-1}}(f^n_{m_n})).
  \]
  This cocycle is invariant under the dual action $\hat{\alpha}$.  Hence we
  obtain another cyclic $(n+1)$-cocycle $i_D \hat{\phi}^\Z$ on $\Z
  \ltimes_\sigma \smooth{A}$, where $D$ is the generator of $\hat{\alpha}$.

  In this setting Corollary~\ref{cor:period-dual-action-pairing-der-cycle-k-group}
  says that the image of $\phi$ under $\partial^*$ pairs with the $K$-group as
  the same way as $i_D\hat{\phi}^\Z$ does.  The strong Morita equivalence~\eqref{eq:mor-equiv-fixed-t-action-from-discr-spec-r-action} becomes the one
  between $\Z \ltimes_\alpha \smooth{A}$ and the crossed product by the
  suspension flow on the mapping torus of $\alpha$.  From this point of view, the
  correspondence of Corollary~\ref{cor:period-dual-action-pairing-der-cycle-k-group} is a generalization of
  the correspondence in Example~\ref{example:pv-nest-corr-torus-transl}.
\end{rmk}

\section{Invariance of the fundamental cocycle}
\label{sec:invar-fund-class}

We turn back to the setting of Section~\ref{sec:preliminaries}, so let $F$ be a foliation on $M$ and $T$ a transverse submanifold of $M$.  Let $q'$
denote the total dimension $q(q+1)/2 + q - 1$ of $\SMet(T)$.  We consider the transverse fundamental
cocycle~\cite{MR866491} of the groupoid $\SMet(G_T)$.  It is defined as the
following cyclic $q'$-cocycle:
\begin{equation}
  \label{eq:defn-transv-fund-cocycle}
  \phi(f^0, \ldots f^{q'}) = \int_{\gamma_0 \cdots \gamma_{q'} \in \SMet T}
  f^0_{\gamma_0} d f^1_{\gamma_1} \cdots d f^{q'}_{\gamma_{q'}}.
\end{equation}

In the following we analyze the pairing of $\hat{\phi}$ with $K_{q'}(\R
\ltimes_{\sigma} C^*\SMet G_T)$ and its relation to the structure of the
von Neumann algebra $W(M; F)$.  Following the method of~\cite{MR866491}, we will
construct an invariant complex measure on the flow of weights when the
cyclic cocycle $\hat{\phi}$ pairs nontrivially with the $K$-group of $\R
\ltimes_{\sigma^{(\omega)}} C^*_r(\SMet(G_T)) \simeq C^*_r(\Met( G_T))$.

\begin{prop}
  Suppose that the foliation $(M; F)$, the transversal $T$, and the transverse density $\omega$ satisfy the condition
  {\condprojinv}.  Then the cocycle $\phi$ of~\eqref{eq:defn-transv-fund-cocycle} is invariant under the $1$-parameter group
  $\sigma_t$.
\end{prop}

\begin{proof}
  The derivative of $\sigma_t^*\phi$ is given by
 \begin{multline*}
 \frac{d \sigma_t^* \phi}{d t}(f^0, \ldots, f^{q'}) = \\
 \sum_{1 \le j \le q'} \int_{\gamma_0 \cdots \gamma_{q'} \in \SMet T} i f^0(\gamma_0) f^j(\gamma_j) d \log(\delta(\gamma_j)) d f^1(\gamma_1) \cdots \widetilde{d f^j(\gamma_j)} \cdots d f^{q'}(\gamma_{q'}) \\
 + \sum_{0 \le j \le q'} \int_{\gamma_0 \cdots \gamma_{q'} \in \SMet T} \log(\delta(\gamma_j)) f^0(\gamma_0) d f^1(\gamma_1) \cdots d f^{q'},
 \end{multline*}
 where $\widetilde{d f^j(\gamma_j)}$ stands for the omission of the corresponding
 term. 

 By assumption \condprojinv, $\delta(\gamma)$ is constant along transverse
 movement.  This implies $d \log(\delta(\gamma_j)) = 0$ for each term in the first part
 of the right hand side.  One also has $\gamma_0 \cdots \gamma_{q'} \in \SMet T$
 implies $\sum_{0 \le j \le q'} \log(\delta(\gamma_j)) = 0$.  This implies that
 the second part is also trivial.
\end{proof}

\begin{rmk}
\label{rmk:godbillon-vey-num-triv}
Suppose that the codimension $q$ is equal to $1$.  The Godbillon--Vey class of $(M; F)$ can be obtained from the derivation of the transverse fundamental class with respect to the action of modular automorphism group~\cite{MR866491}*{Lemma~7.6}.  In particular, the Godbillon--Vey number has to be $0$ under our assumption, but the von Neumann algebra can be of type III as in Example~\ref{exmpl:hyperb-torus}.
\end{rmk}

Connes~\cite{MR866491} showed that the fundamental cocycle of the groupoid $\Met G_T$ extends to a spectral subalgebra of $C^*_r(\Met G_T)$.  We show that a similar property holds for $\SMet G_T$ although it does not have the almost isometric property; that is, we show that $\phi$ extends to a subalgebra of $C^*_r(\SMet G_T)$ which has the same K-groups as $C^*_r(\SMet G_T)$, relying on a deep result on Oka's principle due to J.\mbox{-}B.~Bost~\cite{MR1062964}.

We briefly recall the relevant constructions from~\cite{MR866491}.  Let $E_0$ be a $\SMet G_T$-equivariant real vector bundle over $\SMet T$, endowed with a fiberwise metric.  We do not assume that the metric is invariant under $\SMet G_T$.  Let $E = E_0 \otimes \C$ be its complexification with the Hermitian inner products $(\langle -, - \rangle_{E_x})_{x \in \SMet G_T}$, and $\Gamma_c(\SMet G_T; r^* E)$ be the space of compactly supported continuous maps $f\colon \SMet G_T \rightarrow E$ satisfying $f(\gamma) \in E_{r(\gamma)}$.  Then, $\Gamma_c(\SMet G_T; r^* E)$ admits a right $C_c(\SMet G_T)$-valued inner product
\[
\langle \xi, \eta \rangle_{C^*_r(\SMet G_T)}(\gamma) = \langle \xi(\gamma), \eta(\gamma) \rangle_{E_{r(\gamma)}}.
\]
This allows us to take the completion of $\Gamma_c(\SMet G_T; r^* E)$ as a  Hilbert C$^*$-module over $C^*_r(\SMet G_T)$, denoted by $C^*_r(\SMet G_T; E)$.  The left convolution defines a closable (non-$^*$-)homomorphism
\[
\lambda_E\colon C^*_r(\SMet G_T) \rightarrow \mathcal{L}(\Gamma_c(\SMet G_T; r^* E)),
\]
whose domain is $C_c(\SMet G_T)$.

\begin{lem}\label{lem:eqvr-bdl-proj-inv-met}
Let $E_0$ be an $\SMet(G_T)$-equivariant real vector bundle, and $R$ be a real number.  Suppose that $E_0$ is endowed with a metric satisfying
\begin{equation}\label{eq:eqvr-bdl-proj-inv-met}
\langle \gamma(\xi), \gamma(\eta) \rangle_{r(\gamma)} = \delta(\gamma)^{2 R} \langle \xi, \eta \rangle_{s(\gamma)}, \quad (\gamma \in \SMet(G_T), \xi, \eta \in (E_0)_{s(\gamma)}).
\end{equation}
Then, for any $f \in C_c(\SMet G_T)$, one has $\norm{\lambda_E(f)} = \norm{\sigma^{(\omega)}_{-i R}(f)}_{C^*_r(\SMet G_T)}$.
\end{lem}

\begin{proof}
It is enough to find an operator $\tilde{S}$ on $C_c(\SMet G_T; r^* E)$ satisfying the two conditions
\begin{align*}
\langle \tilde{S}(\xi), \tilde{S}(\xi) \rangle_{C^*_r(\SMet G_T)} &= \langle \xi, \xi \rangle_{C^*_r(\SMet G_T)},&
\tilde{S}(\lambda_E(f) \xi) &= \tilde{S}(\xi) \sigma_{i R}(f^*).
\end{align*}
We define $\tilde{S}$ by $\tilde{S}(\xi)(\gamma) = \delta(\gamma)^{-R} \gamma.\overline{\xi(\gamma^{-1})}$.  Since $\xi(\gamma^{-1}) \in E_{r(\gamma^{-1})} = E_{s(\gamma)}$ and $E$ has a distinguished real subbundle $E_0$, this is a well-defined transformation on $C_c(\SMet G_T; r^* E)$.

On one hand, if $\gamma, \gamma', \gamma'' \in \SMet G_T$ satisfies $\gamma^{-1} = \gamma' \gamma''$, one has
\[
\delta(\gamma)^{-R} \gamma (\gamma' \overline{\xi(\gamma'')}) = \delta(\gamma')^{R} \tilde{S}(\xi)((\gamma'')^{-1}).
\]
On the other hand, by definition of the right convolution, one has
\[
(\eta \sigma_{i R}(f))(\gamma) = \sum_{\bar{\gamma}'' \bar{\gamma}' = \gamma} \delta(\bar{\gamma}')^R \eta(\bar{\gamma}'') f(\bar{\gamma}').
\]
Combining these two, we obtain $\tilde{S}(\lambda_E(f) \xi) = \tilde{S}(\xi) \sigma_{i R}(f^*)$.  Next, using the condition~\eqref{eq:eqvr-bdl-proj-inv-met}, one computes
\[
\langle \tilde{S}(\xi), \tilde{S}(\xi) \rangle_{C^*_r(\SMet G_T)}(\gamma) = \langle \xi(\gamma^{-1}), \xi(\gamma^{-1}) \rangle = \langle \xi, \xi \rangle_{C^*_r(\SMet G_T)}(\gamma).
\]
This completes the proof.
\end{proof}

\begin{rmk}\label{rmk:eqvr-tang-sub-quot-decomp}
Let $K = \Tang \ker p$ be the tangent bundle of the fibers for the fiber bundle $p\colon \SMet T \rightarrow T$.  Since $\GL_k(\R) / O_k(\R)$ has an $\GL_k(\R)$-invariant metric, $K$ is a $\SMet G_T$-equivariant real vector bundle with an invariant metric.  We define $D_0' = p^* \Tang T \oplus K$, and $D_0 = \oplus_{k = 1}^{q'} (D_0')^{\otimes k}$.  The tautological metric on $p^* \Tang T$ satisfies the assumption of Lemma~\ref{lem:eqvr-bdl-proj-inv-met} with $R = 1$.  It follows that $E'_0$ is a direct sum of vector bundles satisfying this assumption with $R = 0, 1, \ldots, q'$. 
\end{rmk}

Let $E_0'$ be the tangent bundle of $\SMet T$, endowed with a natural action of $\SMet G_T$.  We fix a subbundle $N \subset E_0'$ which satisfies $K \cap N = 0$ and $K + N = E_0'$.  This allows us to define a metric on $E_0$ which is uniquely defined by the conditions $N \perp K$, $N \simeq p^* \Tang T$.  In other words, $D_0'$ and $E_0'$ are the same as metric bundles over $\SMet T$, only $\SMet G_T$ acts differently.  As in the case of $D_0$, we put $E_0 = \oplus_{k = 1}^{q'} (E_0')^{\otimes k}$.

\begin{prop}\label{prop:subalg-transv-fund-def}
There is a Banach subalgebra $\mathcal{A}$ of $C^*_r(\SMet G_T)$ such that
\begin{enumerate}
\item $\mathcal{A}$ is contained in the domain of the closure of $\lambda_E$,
\item the inclusion $\mathcal{A} \rightarrow C^*_r(\SMet G_T)$ induces isomorphisms of the K-groups,
\item the action $\sigma$ restricts to a strongly continuous action on $\smooth{A}$.
\end{enumerate}
\end{prop}

\begin{proof}
We freely use the notation of Remark~\ref{rmk:eqvr-tang-sub-quot-decomp} in this proof.

First, put $\bnorm{f}_0 = \sup_{\absolute{R} \le q'} \norm{ \sigma_{i R}(f)}_{C^*_r(\SMet G_T)}$.  By Lemma~\ref{lem:eqvr-bdl-proj-inv-met}, $\lambda_{D}$ is continuous with respect to $\bnorm{-}_0$.  Moreover,~\cite{MR1062964}*{Theorem~1.1} implies that the closure $\mathcal{A}_0$ of $C_c(\SMet G_T)$ by $\norm{-}_0$ has the same K-groups as $C^*_r(\SMet G_T)$.  Meanwhile, by Lemma~\ref{lem:eqvr-bdl-proj-inv-met} and Remark~\ref{rmk:eqvr-tang-sub-quot-decomp}, $\lambda_D$ defines a contractible representation of $\mathcal{A}_0$ on $C^*_r(\SMet G_T; D)$.

Next, $K$ is a $\SMet G_T$-invariant subbundle of $E_0'$.  Thus, the action of $G_T$ on $E_0'$ can be represented as a triangular block matrix with orthogonal matrices in the block diagonals, with respect to the decomposition $E_0' \simeq K \oplus p^* \Tang T$~\cite{MR866491}*{Lemma~5.2}.  This implies that there is a linear map $P$ from $C_c(\SMet G_T)$ to the space of linear transformations on $\Gamma_c(\SMet G_T; E')$ such that $\lambda_{E'}(f) = \lambda_{D'}(f) + P(f)$.  As in the proof of \cite{MR866491}*{Lemma~3.4}, $\lambda_{E'}(f)$ is conjugate to $\lambda_{D'}(f) + \epsilon P(f)$ for any $\epsilon > 0$.  Then, we can find linear maps $P_1, \ldots, P_{q'}$ from $C_c(\SMet G_T)$ to the space of linear transformations on $\Gamma_c(\SMet G_T; E)$ such that $\lambda_E(f)$ is conjugate to $\lambda_\epsilon(f) = \lambda_D(f) + \sum_{k = 1}^{q'} \epsilon^k P_k(f)$ for any $\epsilon > 0$.

Now, we can imitate the proof of~\cite{MR866491}*{Proposition~3.5}.  We define $\mathcal{A}$ to be the intersection of $\mathcal{A}_0$ and the domain of the closure of $\lambda_E$.  It remains to show that $\mathcal{A}$ is closed under holomorphic functional calculus in $\mathcal{A}_0$.  For this, it is enough to show that, whenever $f \in \mathcal{A}$ satisfies $\bnorm{f}_0 < 1$, the inverse of $1 + f$ lies in $\mathcal{A}_0$~\cite{MR866491}*{proof of Lemma~1.2}.  Suppose that this assumption holds.  Then, $\norm{\lambda_D(f)}$ is also smaller than $1$.  Since the mapping $(0, 1) \ni \epsilon \mapsto \lambda_\epsilon(f)$ is continuous and has the limit $\lambda_D(f)$ at $\epsilon = 0$, one has to have $\norm{\lambda_\epsilon(f)} < 1$ at certain $\epsilon$.  Now, the algebra $\smooth{A}$ is equal to the closure of $C_c(\SMet G_T)$ with respect to the norm $\bnorm{f} = \sup(\bnorm{f}_0, \norm{\lambda_\epsilon(f)})$ for this $\epsilon$.  Thus, the inverse of $1 + f$ exists in $\mathcal{A}$.

Finally, it remains to verify that $\sigma$ restricts to a strongly continuous action on $\smooth{A}$.  By definition of $\bnorm{-}_0$, $\sigma$ restricts to a continuous action on $\smooth{A}_0$.  Similarly, the multiplier $\delta(\gamma)$ defines an strongly continuous one-parameter automorphism group on $\mathcal{L}(C^*_r(\SMet G_T; E))$, which extends $\sigma$ via $\lambda_E$.  Thus, we obtain that $\sigma$ is strongly continuous on $\mathcal{A}$.  This completes the proof.
\end{proof}

We are ready to apply the results of Section~\ref{sec:peri-acti-cycl} to the transverse fundamental cocycle $\phi$.  By~\cite{MR866491}*{Theorem~3.7}, $\phi$ is an $q'$-trace with respect to the norm $\bnorm{f}$ for $\epsilon = 1$ in the above proof.  Thus, there exists a Fr\'{e}chet subalgebra $\smooth{A}^\infty_0$ of $\smooth{A}$ which is closed under holomorphic functional calculus, and to which $\phi$ extends as a $q'$-cyclic cocycle~\cite{MR866491}*{Section~2}.  The symmetry of the seminorms used to define $\smooth{A}^\infty_0$ implies that $\sigma$ acts strongly continuously on $\smooth{A}^\infty_0$.  We define $\smooth{A}^\infty$ to be the subalgebra of $\smooth{A}^\infty_0$ consisting of the smooth elements of $\sigma$.  This is again stable under holomorphic functional calculus in $\smooth{A}^\infty_0$.  Since we made sure the same K-groups when we pass to subalgebras from $C^*_r(\SMet G_T)$ to $\smooth{A}^\infty$ at each step, we obtain that the inclusion $\smooth{A}^\infty \rightarrow C^*_r(\SMet G_T)$ induces the isomorphisms of K-groups.

\subsection{Anabelian cocycles}
\label{sec:anabelian-cocycles}

\begin{dfn}[\citelist{\cite{MR866491}*{Theorem~7.14}\cite{MR1303779}*{Section~3.6.$\gamma$ Definition~18}}]
  Let $A$ be a C$^*$-algebra, $\smooth{A}$ be its dense algebra, and $B$ be a unital subalgebra of the center $Z(\Mult(A))$ of the multiplier
  algebra of $A$.  A cyclic $k$-cocycle $\psi$ on $B  \smooth{A}$ is said to be
  $B$-\textit{anabelian} if the $(k + 1)$-linear map on $\smooth{A}$
  \[
  \psi_g (f^0, \ldots, f^k) = \psi(g f^0, f^1, \ldots, f^k)
  \]
  is a cyclic $k$-cocycle for any $g \in B$.
\end{dfn}

\begin{lem}
  \label{lem:anabelian-van-first-var-cond}
  Let $A$, $\smooth{A}$, $B$, be as above, and let $\psi$ be a $k$-cocycle on $\smooth{A}$.  Suppose that $\psi$ extends to a cyclic cocycle on
  $(\smooth{A} + B)^{k + 1}$ and satisfies
 \begin{equation}
   \label{eq:anabelian-first-var-van-cond}
  \psi(g, f^1,\ldots, f^k) = 0
  \end{equation}
  for any $g \in B$.  Then $\psi$ is $B$-anabelian.
\end{lem}

\begin{proof}
  The Hochschild cocycle property of $\psi_g$ follows from that of $\psi$ and
  from the fact that $g$ commutes with any element of $\smooth{A}$.  Hence it
  remains to show that $\psi_g$ satisfies the cyclicity condition.

  The cyclicity of $\psi$ and~\eqref{eq:anabelian-first-var-van-cond}
  implies
  \begin{equation}
    \label{eq:anabelian-any-place-van-cond}
  \psi(f^j, \ldots, f^k, g, f^1, \ldots, f^{j - 1}) = 0
  \end{equation}
  for any $1 \le j \le k$.  The Hochschild cocycle condition implies
  \begin{equation*}
  \begin{split}
  0 &= b\psi(f^0, g, f^1, \ldots, f^k) \\
  &= \psi(f^0 g, f^1, \ldots, f^k) - \psi(f^0, g f^1, f^2, \ldots, f^k) \\
  &+ \sum_{j=2}^k (-1)^j \psi(f^0, g, f^1, \ldots, f^{j-1} f^j,\ldots, f^k) + (-1)^k \psi(f^k f^0, g, f^1, \ldots, f^{k-1}).
  \end{split}
  \end{equation*}
  Combined with~\eqref{eq:anabelian-any-place-van-cond} and the cyclicity
  condition of $\psi$, this implies the cyclicity condition of $\psi_g$.
\end{proof}

\begin{prop}
  Suppose that $(M; F)$, $T$, and $\omega$ satisfy the conditions {\condprojinv}
  and {\conduniterg}.  Then the dual $\hat{\phi}$ on $\R \ltimes_{\sigma}
  C^\infty_c(\SMet G_T)$ of the transverse fundamental cocycle is $C^\infty(\Det
  T)^{G_T}$-anabelian.
\end{prop}

\begin{proof}
  By Lemma~\ref{lem:anabelian-van-first-var-cond}, it is enough to show the
  equality
  \[
  \hat{\phi}(g, f^1, \ldots, f^{q'}) = 0
  \]
  for any $G_T$-invariant function $g$ on $\Det T$.  Recall the
  decomposition $\Det(T) \simeq T \times \R_{>0}$ determined by the choice of
  $\omega$.  Then we can write
  \begin{equation}
 \label{eq:int-rep-psi-g-fs}
  \psi(g, f^1,\ldots,f^{q'}) = \int_{\R_{>0}} \frac{ d t}{t} \int_{\SMet(T) \times \ensemble{t}} g(p(x), t) d f^1 \cdots d f^{q'},
  \end{equation}
  where $p\colon \Det(T) \rightarrow T$ is the natural projection.

  By the condition {\conduniterg}, the function $g(y, t)$ on $T \times \R_{>0}
  \simeq \Det(T)$ is a constant function for each $t$.  Hence one has
  \[
  \int_{\SMet(T)} g(p(x), t) d f^1 \cdots d f^{q'} = 0,
  \]
  which proves that the integral in the right hand side of~\eqref{eq:int-rep-psi-g-fs} is trivial.
\end{proof}

\begin{thm}
  \label{thm:inv-meas-flow-weights}
  Suppose that the $(M; F)$, $T$, and $\omega$ satisfy the conditions
  {\condprojinv} and {\conduniterg}.  For each $x$ in $K_{q'}(C^*_r(\Met G_T))$,
  there is a $\theta$-invariant normal functional $\nu_x$ on the center of $\R
  \ltimes_\sigma W(M; F)$ satisfying $\nu_x(1) = \pairing{\hat{\phi}}{x}$.
\end{thm}

\begin{proof}
  When $f \in C^\infty(\Det(T))^{G_T}$, we define
 \[
 \nu_x(f) = \pairing{ \hat{\phi}_f }{x}.
 \]
 The right hand side is normal as a functional defined on a weakly dense
 subspace of $L^\infty(\Met(T))^{G_T}$.  Moreover it satisfies $\nu_x(1) =
 \pairing{\hat{\phi}}{x}$.  Hence it extends to a normal functional over
 $L^\infty(\Met(T))^{G_T}$.

 By $\hat{\sigma}|_{\R \ltimes C^*_r(G_T)} = \theta$ and the
 $\hat{\sigma}$-invariance of the cocycle $\hat{\phi}$, we have
 \[
 \pairing{\theta_t^* \hat{\phi}}{x} = \pairing{\hat{\phi}}{\hat{\sigma}_{t}(x)}.
 \]
 Since $\hat{\sigma}_t(x)$ is a constant family in $K_*(C^*\Met(G_T))$, one
 has $\pairing{\hat{\phi}}{\hat{\sigma}_{t}(x)} = \pairing{\hat{\phi}}{x}$.
 This proves the invariance of $\nu_x$ under $\theta$.
\end{proof}

\begin{cor}
With the same assumption as in Theorem~\ref{thm:inv-meas-flow-weights}, assume moreover that there exists an element $x \in K_{q'}(C^*_r(\Met G_T))$ such that $\pairing{\hat{\phi}}{x} \neq 0$.  Then $W(M; F)$ is a factor of type \textnormal{III}.
\end{cor}

\begin{proof}
The flow of weights of a semifinite factor is isomorphic to the translation of $\R$ on $L^\infty(\R)$.  Hence there cannot be an invariant normal functional on $Z(\R \ltimes_\sigma W(M; F))$ in that case~\cite{MR866491}.
\end{proof}

\section{Period of flow and \texorpdfstring{$K$}{K}-cycles}
\label{sec:per-flow-k-cycle}

In this section we apply the consideration of Section~\ref{sec:peri-acti-cycl} to the situation of
Theorem~\ref{thm:inv-meas-flow-weights} and calculate the possible values of
pairing between the dual fundamental cocycle and the $K$-group when the von
Neumann algebra $W(M; F)$ is a factor of type III$_\lambda$ for some $0 < \lambda < 1$.  We will use Corollary~\ref{cor:saturated-action-pimsner-voiculescu-nest} to relate the $K$-theory pairing of the dual of the fundamental cocycle with $K_{q'}(C^*_r(\Met(G_T)))$ to that
of the fundamental cocycle with $K_{q'}(C^*_r(\SMet(G_T^{(u)})))$.

In order to construct nontrivial elements in the $K$-group of groupoid algebras,
we consider the group $K^{\topol}_*(G_X)$ and the assembly
map~\citelist{\cite{MR1769535}\cite{MR866491}\cite{MR679730}}
\[
\mu \colon K^{\topol}_*(G_X) \rightarrow K_*(C^*_r(G_X))
\]
for the base spaces $X = \SMet(T)$ and $\Met(T)$.

A cycle $c$ in $K^{\topol}_0(G_X)$ is represented by a quadruple $(N, f, E, D)$: $N$ is a manifold endowed with a map $f\colon N \rightarrow X$\ and a proper action of $G_X$ with respect to $f$.  Furthermore $E$ is a $\Z_2$-graded $G_X$-equivariant vector bundle on $N$, and $D = (D_x)_{x \in X}$ is a family of odd elliptic operators with coefficient $E$ on the fibers of $f$.  Here, $D$ is assumed to satisfy the $G_X$-equivariance condition $\gamma D_{s\gamma} \gamma^{-1} = D_{r\gamma}$ for
$\gamma \in G$.  Given such data $c$, the $G$-index
$\ind_G(D_*)$ of the family $(D_x)_{x\in X}$ defines an element $\mu(c)$
in $K_0(C^*_r(G_X))$.

Similarly, the pairs of $G_X$-equivariant map $f\colon N \rightarrow X$ and
equivariant longitudinal elliptic operators over $N \times \R$ define odd cycles which are elements of $K^{\topol}_1(G_X)$, and elements in $K_1(C^*_r(G_X))$ via the assembly map.

Now, let us consider the case of the base space $X = \SMet(T)$.  In the following we shall describe the natural map
\[
\Phi^{\topol}\colon K^{\topol}_0(G_X) \rightarrow K^{\topol}_1(\Det(G_X)),
\]
which is compatible with the Connes--Thom isomorphism
$\Phi$ from $K_0(C^*_r(G_X))$ to $K_1(\R \ltimes_{\sigma^{(\omega)}} C^*_r(G_X))$ via the assembly map $\mu$.

Suppose that an element $x$ of $K^{\topol}(G_X)$ is represented by $(N, f, E, D_*)$ of
a $G_X$-equivariant map $f\colon N \rightarrow X$ and a family $(D_x)_{x\in X}$
of odd elliptic operators on a $G_X$-equivariant graded vector bundle $E$.  Then $N \times \R$ admits an action
of $G$ defined by
\[
\gamma (x, r) = (\gamma x, \delta^{(\omega)}(\gamma) t).
\]
The induced map
\[
\tilde{f}\colon N \times \R \rightarrow X \times \R \simeq \Det X \simeq
\Met(T)
\]
is $G_X$-equivariant under this action.

Put $Y = X \times \R^2$ and let $\tilde{\smooth{E}}$ be the $G_X$-equivariant
vector bundle $\tilde{\smooth{E}}_{(x, v)} = \smooth{E}_x \otimes \C^2$.  For
each $(x, v) \in X \times \R^2$, put $\tilde{D}_{(x, v)} = D_x \otimes 1 +
\epsilon \otimes D_{\R^2}$, where $\epsilon$ is the grading operator on
$\smooth{E}$ and $D_{\R^2}$ is the Dirac operator on $\R^2$.  Then it defines a $G_X$-elliptic operator and the data $(Y,
\tilde{f}, (\tilde{D}_{(x, v) \in X \times \R^2}))$ is a family of equivariant
fiberwise elliptic operators parametrized by $Y = \Det(X) \times \R$.  Let $\Phi^{\topol}(x)$ denote this element of $K^{\topol}_1(\Det(G_X))$.

Similarly, we can define $\Phi^{\topol}\colon K^{\topol}_1(G_X) \rightarrow K^{\topol}_0(G_X)$ by the Bott periodicity.

\begin{lem}
  \label{lem:topol-desc-connes-thom}
  We have the equality
  \[
  \mu(\Phi^{\topol}(N, f, E, D_*)) = \Phi(\mu(N, f, E, D_*))
  \]
  in $K_{*+1}(\R \ltimes_{\sigma^{(\omega)}} C^*_r(G_X))$ for any cycle $(N,
  f, E, D_*)$ of $K^{\topol}_*(G_X)$.
\end{lem}

\begin{proof}
Let $\alpha^{(t)}$ denote the action $\alpha^{(t)}_s(x) = \sigma^{(\omega)}_{ts}(x)$ of $\R$ over $C^*_r(G_X)$, and $A = (A_t)_{t \in [0, 1]}$ be the continuous field of C$^*$-algebras over $[0, 1]$ whose fiber $A_t$ at $t \in [0, 1]$ is given by the crossed product $\R \ltimes_{\alpha^{(t)}} C^*_r(G_X)$.  Then the evaluation map $e_0\colon A \rightarrow A_0$ at $t = 0$ is an \KK-isomorphism and the one $e_1\colon A \rightarrow A_1$ at $t = 1$ can be considered as a \KK-morphism
  \[
  \Lambda\colon C_0(\R) \otimes C^*_r(G_X) \simeq A_0 \rightarrow A_1.
  \]
The Connes--Thom isomorphism $K_*(C^*_r(G_X)) \rightarrow K_{*+1}(\R \ltimes_{\sigma^{(\omega)}} C^*_r(G_X))$ is given by the composition of the Bott periodicity isomorphism
\[
 K_*(C^*_r(G_X)) \rightarrow K_{*+1}( C_0(\R) \otimes C^*_r(G_X)) \simeq K_{*+1}(A_0)
\]
and the map induced by $\Lambda$.

Put $Z = X \times \R^2 \times [0, 1]$ and consider the action of $G_X$ on $Z$ by
\[
 \gamma.(x, v, t) = (\gamma.x, (t \log(\delta^{(\omega)}) + v_1, v_2), t).
\]
Then the groupoid algebra $C^*_r(G_X \ltimes Z)$ is isomorphic to $A$ as a $C([0, 1])$-algebra.  Given a geometric cycle $z = (N, f, E, D_*)$ for $G_X$, let $\tilde{z}$ be the geometric cycle $(N \times \R^2 \times [0, 1], p_1^*(D_*))$ over $G_X \ltimes Z$.  Then one has $e_0(\mu(\tilde{z})) = \mu(z)$ and $e_1(\mu(\tilde{z})) = \Phi(\mu(z))$, which proves the assertion.
\end{proof}

Let $c = (N, f, E, D)$ be a cycle in $K^{\topol}_*(\Det(G_X))$.  We take the map $N \rightarrow \R$ given by the composition of $f$ and the projection $\Det X \rightarrow \R$ induced by $\omega$.  Then, the
inverse image $N_0$ of $X \times \ensemble{0}$ in $N$ is a submanifold of
codimension $1$, and it has a $G_X^{(u)}$-equivariant map into $X$.  Moreover, $E$ and $D$ restricts to $N_0$.  The data $\iota^*(c) = (N_0, f|_{N_0}, E|_{N_0}, D|_{N_0})$ defines an element of $K^{\topol}_*(G_X)$.  This correspondence defines a well-defined map $\Psi^{\topol}\colon K^{\topol}_*(\Det(G_X)) \rightarrow K^{\topol}_*(G_X)$.

It can be easily seen that this map corresponds to the homomorphism $\Psi$ of~\eqref{eq:mapping-torus-periodization-defn} via $\mu$.

\begin{lem}
  \label{lem:topol-pim-voi-boundary-map}
  We have the equality
  \[
  \mu(\Psi^{\topol}(N, x)) = \Psi(\mu(N, x))
  \]
  in $K^{\topol}(G_X^{(u)})$ for any $(N, x) \in K^{\topol}(\Det(G_X))$.
\end{lem}

Let $A$ be the C$^*$-algebra $C^*_r(\SMet G_T)$.  As a consequence of Lemmas~\ref{lem:topol-desc-connes-thom} and~\ref{lem:topol-pim-voi-boundary-map}, we have the following commutative
  diagram:
 \[
 \xymatrix{
 K^{\topol}_*\big(\SMet G_T\big) \ar[r]^-{\Phi^{\topol}} \ar[d]_{\mu} &
 K^{\topol}_{*+1}(\Met(G_T)) \ar[r]^-{\Psi^{\topol}} \ar[d]_{\mu} &
 K^{\topol}_{*+1}\big(\SMet G_T^{(u)}\big)  \ar[d]_{\mu} \\
 K_*(A) \ar[r]^-{\Phi} \ar@/_1em/[rr]_{\partial} & K_{*+1}(\R
 \ltimes_{\sigma^{(\omega)}} A) \ar[r]^-{\Psi} & K_{*+1}\big(C^*_r\big(\SMet G_T^{(u)}\big)\big).
 }
 \]

\begin{thm}
  \label{thm:iii-lambda-proj-inv-per-calc-trans-fund}
  Let $(M; F)$ be a foliation, $T$ a transversal and $\omega$ a smooth density
  on $T$, satisfying the conditions {\condprojinv} and {\conduniterg}.  Suppose
  that the von Neumann algebra $W(M; F)$ is of type \textnormal{III}$_\lambda$ for $0 <
  \lambda < 1$.  Then we have
 \[
 \ensemble{\pairing{i_D\phi}{\mu(x)} \suchthat x \in K^{\topol}(\SMet(G_T)) } \subset \log(\lambda) \Q.
 \]
 We also have
 \[
  \log(\lambda) \Z \subset \ensemble{ \pairing{i_D\phi}{y} \suchthat y \in K_{q+1}(C^*_r(\SMet(G_T)))}.
 \]
\end{thm}

\begin{proof}
  By Lemma~\ref{lem:proj-inv-density-s-set-calc} and the assumption on the type
  of $W(M; F)$, the spectrum of $\sigma^{(\omega)}$ agrees with $\log(\lambda)
  \Z$.  By {\conduniterg}, the action $\sigma^{(\omega)}$ is saturated.
  
 By Proposition~\ref{prop:subalg-transv-fund-def} and the following remark, we can find a $\sigma$-invariant Fr\'{e}chet subalgebra $\smooth{A}^\infty$ of $C^*_r(\SMet G_T)$ with the same K-groups, such that $\phi$ extends.  Hence we may apply Corollary~\ref{cor:saturated-action-pimsner-voiculescu-nest}, and obtain the equality
 \begin{equation}
 \label{eq:iii-lambda-case-dual-fund-pair}
 \pairing{i_D\phi}{y} = \log(\lambda) \pairing{\phi|_{\overline{C^\infty_c(\Met G_T^{(u)})}}}{\partial(y)}
 \end{equation}
 for any $y \in K_*(C^*_r(G_T))$.

 For any $(N, f, E, D) \in K^{\topol}_q(G_T^{(u)})$,
 we have
 \[
 \pairing{\phi}{\mu(N, z)} = \int_N \Todd(f) \ch(\sigma_D) \in \Q,
 \]
 where $\ch(\sigma_D)$ denotes the Chern character of the symbol of $D$ in the compact support cohomology group $H^*_c(N, \C)$, and $\Todd(f)$ is the relative Todd class of $f$.  This and~\eqref{eq:iii-lambda-case-dual-fund-pair} imply the first assertion.

 Let $[G_T^{(u)}]^* \in K_q(C^*_r(G_T^{(u)}))$ be the dual fundamental cocycle
 defined in~\cite{MR679730}*{Section~8}.  It is represented by $\mu(D^q, x)$,
 where $D^q \rightarrow T$ is an immersion of the $q$-dimensional open disk and
 $x \in K^q(D^q)$ is the generator of $K^q(D^q) \simeq \Z$.

 We show that $[G_T^{(u)}]^*$ is in the image of $K_{q+1}(C^*_r(G_T))$, which is
 equivalent to
 \begin{equation}
 \label{eq:trans-fund-dual-class-inv-under-dual-action}
 (1-\hat{\sigma})(\Xi([G_T^{(u)}]^*)) = 0
 \end{equation}
 for the isomorphism $\Xi\colon K_*(C^*_r(G_T^{(u)})) \rightarrow K_*(\T
 \ltimes_{\sigma^{(\omega)}} C^*_r(G_T))$ induced by the strong Morita
 equivalence.

 Let $G_T \times_\omega \Z$ be the groupoid whose object set is $T \times \Z$,
 arrow set is $G_T \times \Z$, and structure maps are given by
 \begin{align*}
 s(\gamma, n) &= (s\gamma, n),& r(\gamma, n) &= \left(r\gamma, n + \frac{\log (\delta^{(\omega)}(\gamma)}{\log \lambda}\right).
 \end{align*}
 Then $G_T \times_\omega \Z$ is strongly Morita equivalent to $G_T^{(u)}$ as a
 groupoid, and one has an isomorphism
 \[
 C^*_r(G_T \times \Z) \simeq \T \ltimes_{\sigma^{(\omega)}} C^*_r(G_T)
 \]
 of algebras.  Under this isomorphism the dual action $\hat{\sigma}$ on $\T
 \ltimes_{\sigma^{(\omega)}} C^*_r(G_T)$ corresponds to the action of $\Z$ on $G_T
 \times_\omega \Z$ given by $\alpha_m(\gamma, n) = (\gamma, n + m)$.

 If $\Xi([G_T^{(u)}]^*)$ is represented by a immersion $D^q \rightarrow T \times
 \ensemble{n}$ for some $\iota\colon D^q \rightarrow T$ and $n \in \Z$, then $\alpha(\Xi([G_T^{(u)}]^*))$ is represented by $(\iota, n+1)\colon D^q
 \rightarrow T \times \ensemble{n+1}$.  By replacing $D^q$ with a smaller disk
 if necessary, we may assume that the image of $\iota$ is contained in a domain
 of an element $\gamma \in G_T$ satisfying $\delta^{(\omega)}(\gamma) =
 \lambda$.  Then $\gamma\circ\iota$ is an immersion of $D^q$ into $T$ and one has the equality
 \begin{equation}
 \label{eq:fund-class-dual-action-image}
 (\gamma\circ\iota, 1)_* [D^q]^* = ((\gamma, 0)\circ(\iota, 0))_* [D^q]^* = \gamma (\iota_*([D^q]^*)) \gamma^* = \iota_*[D^q]^*
 \end{equation}
 in $K_q(C^*_r(G_T \times_\omega \Z))$.  By the condition {\conduniterg}, there is
 an element $\gamma_1$ in $G_T^{(u)}$ whose domain contains the image of $\gamma
 \circ \iota$ and codomain lies in the same connected component as the image of
 $\iota$.  With this $\gamma_1$ one has
 \begin{equation}
 \label{eq:fund-class-conn-comp-jump}
 (\gamma\circ\iota, 1)_* [D^q]^* = \gamma_1 ((\gamma\circ\iota, 1))_* [D^q]^*  \gamma_1^* = (\iota, 1)_*[D^q] = \alpha((\iota, 0)_*[D^q]^*).
 \end{equation}

 Combining~\eqref{eq:fund-class-dual-action-image} and~\eqref{eq:fund-class-conn-comp-jump}, we obtain $\alpha(\iota_*[D^q]^*) =
 \iota_*[D^q]^*$, which proves~\eqref{eq:trans-fund-dual-class-inv-under-dual-action}.  By $\pairing{\phi}{[G_T^{(u)}]^*} = 1$, \eqref{eq:iii-lambda-case-dual-fund-pair}, and Corollary~\ref{cor:saturated-action-pimsner-voiculescu-nest}, we obtain
 \[
 \log(\lambda) \in \pairing{i_D\phi}{K_{q+1}\big(C^*_r\big(\SMet(G_T)\big)\big)}.
 \]
  This proves the latter half of the assertion.
\end{proof}

\subsection{The \texorpdfstring{$1$}{1}-form corresponding to the dual of fundamental cocycle}
\label{sec:form-cor-dual-fund}

In this last section, we consider the case where $M$ is compact, there is a covering $\pi\colon M' \rightarrow M$ and a holonomy invariant transverse density $\omega$ of $(M', \pi^*F)$, such that $\omega$ is projectively invariant under the deck transformation group of $M'$.  Under this assumption, if we take a transversal $T$ which admits a section $f\colon T \rightarrow M'$, the pullback of $\omega$ on $T$ by $f$ defines a projectively invariant transverse density for $F$.

The Radon--Nikodym cocycle of the deck transformation group with respect to the density $\omega$ defines an $\R_+$-valued group $1$-cocycle.  By taking the natural logarithm and considering the double complex of $\Omega^*(M')$-valued $\Gamma$-cochains, this group cocycle can be regarded as a class $[\log \delta^{(\omega)}]$ of $H^1(M; \R)$.

We also can define $[\log \delta^{(\omega)}]$ by a \v{C}ech $1$-cocycle in the following way.  Let $(U_i)_{i \in I}$ be a covering of $M$ by foliation charts admitting sections $f_i \colon U_i \rightarrow M'$.  Then, for each pair $(i, j)$ of indices there exists a unique element $g_{i j}$ of the deck transformation group.  The scalar $c_{i j}$ characterized by $c_{i j} \omega = g_{i j}^* \omega$ satisfies a \v{C}ech　$1$-cocycle identity with respect to the group law of $\R_+^\times$.  The $\R$-valued \v{C}ech $1$-cocycle $(\log(c_{i j}))_{i j}$ for $(U_i)_{i \in I}$ defines the desired class $[\log \delta^{(\omega)}]$.

For each open set $U$ of $M$, let $\smooth{A}_U$ be the convolution algebra of compactly supported smooth
functions over the full holonomy groupoid of $F|_U$.  On the one hand, any family
$(D_l)_{l \in M/F}$ of leafwise elliptic operators defines an element $\ind(D)$
of the $K_0$-group of $\smooth{A}_M$.  On the other hand, by the strong
Morita equivalence, there is a natural isomorphism $\HC^{n}(C^\infty_c(G_T))$
and $\HC^n(\smooth{A}_M)$ for each $n$.

Next, recall that there is a map of $\Z_2$-graded vector spaces~\cite{MR1303779}*{Section~3.7.$\gamma$}
 \begin{equation}
 \label{eq:chern-char-from-period-cyclic-fol-alg-to-usual-cohom}
 \Lambda\colon \HP^*(\smooth{A}_M) \rightarrow \big(\bigoplus_{k \in q + 2 \Z} H^k(M) \big) \oplus \big(\bigoplus_{k \in 1 + q + 2 \Z} H^k(M) \big)
 \end{equation}
which satisfies
 \[
 \pairing{\ind(D)}{\psi} = \int_M \Todd(\tau \otimes \C)^{-1} \Phi_{H}(\ch([\sigma(D)])) \Lambda(\psi),
 \]
 where $\Phi_{H}$ is the Thom isomorphism $H^*_c(F^*) \rightarrow H^{*-p}_c(M)$
 in cohomology.  Let us recall the construction of $\Lambda$.
 
 We take an open covering $\mathcal{U}
 = (U_i)_{i \in I}$ of $M$.  For each open set $U \subset M$ and $k \in \N$, put
 \[
 \Omega_F^k(U) = (\smooth{A}_U^{\otimes k} \oplus
 \smooth{A}_U^{\otimes k + 1})'.
 \]
 Next, for such $U$ and $k$, let $\Omega^\tau_k(U)$ denote the space of
 $F|_U$-holonomy invariant transverse $k$-currents on $U$.  By choosing a
 transverse bundle $H \subset \Tang M$, one can define a map of complexes
 \begin{equation}
   \label{eq:cyclic-transv-densty-straight-set}
   \Omega^\tau_k(U) \rightarrow \Omega_F^k(U)
 \end{equation}
for any $k$, natural in $U$, and does not depend on the choice of $H$ where one
passes to the cohomology.

Note that when $U_1 \subset U_2$ are open sets of $M$, we have natural maps
$\Omega_F^k(U_2) \rightarrow \Omega_F^k(U_1)$ and $\Omega^\tau_k(U_2)
\rightarrow \Omega^\tau_k(U_1)$.

We then consider the following
 triple complexes $\Gamma_{\mathcal{U}}^{a, b, c}$ and
 ${\Gamma'}_{\mathcal{U}}^{a, b, c}$.
 The first one is given by
 \[
 \Gamma_{\mathcal{U}}^{a, b, c} = \prod_{i_1 < \cdots < i_c}
 \Omega_F^{b - a}(U_{i_1} \cap \cdots \cap U_{i_c}),
 \]
together with the differentials $B$ (which increases the index $a$), $b$ (resp. the index $b$), and the \v{C}ech coboundary map (resp. the index $c$).  Similarly, the second triple complex is given by
 \[
 {\Gamma'}_{\mathcal{U}}^{a, b, c} = \prod_{i_1 < \cdots < i_c}
 \Omega^\tau_{b - a}(U_{i_1} \cap \cdots \cap U_{i_c}).
 \]
 with the de Rham differential $d$ (which increases the index $a$), the zero differential (resp. the index $b$), and the \v{C}ech coboundary map (resp. the index $c$).  Then the map~\eqref{eq:cyclic-transv-densty-straight-set} induces a map of complexes
 \begin{equation}
   \label{eq:cyclic-triple-to-current-triple}
   \Lambda_0 \colon {\Gamma'}_{\mathcal{U}}^{a, b, c} \rightarrow \Gamma_{\mathcal{U}}^{a, b, c},
 \end{equation}
 which is a quasi-isomorphism when each $U_i$ is contained in a foliation chart.

 When $U$ is a foliation chart of
 $F$ homeomorphic to $T \times V$ for $T \simeq \R^q$ and $V \simeq \R^p$, we have the isomorphism $\smooth{A}_U \simeq C^\infty_c(T) \ptensor
 \smooth{K}^\infty$, which implies
 \[
 \HP^n(\smooth{A}_U) \simeq
 \begin{cases}
   \C & (q \equiv n \mod 2)\\
   0 &  (q \not\equiv n \mod 2).
 \end{cases}
 \]

 On the other hand we also have $\Omega^\tau_k(U) = \Omega_k(T)$ for such $U$.  It follows that the cohomology of the complex $\Omega^\tau_*(U)$ is trivial except for the degree $q$, and that $H_q(\Omega^\tau_*(U)) = \C$.  When $(U_i)_{i \in I}$ is a good cover of $M$ by foliation charts, the cohomology of the triple complex ${\Gamma'}_{\mathcal{U}}^{a, b, c}$ is equal to the \v{C}ech cohomology of the orientation sheaf of $\tau$ (up to a shift of degree), which is isomorphic to the usual cohomology group of $M$.

 Hence the cohomology of the triple complex $\Gamma_{\mathcal{U}}^{a, b, c}$ is
 equal to the right hand side of~\eqref{eq:chern-char-from-period-cyclic-fol-alg-to-usual-cohom}.  The restriction homomorphism $\Omega_F^{b - a}(M)
 \rightarrow \prod_{i \in I} \Omega_F^{b - a}(U_i)$ gives the required map $\Lambda$.
 
 Although the manifolds $\Det(M)$, $\Met(M)$ are not compact, the above construction of $\Lambda$ still makes sense for the foliations induced by $F$ on these manifolds.  For example, when $U$ is a foliation chart of $(M; F)$, its inverse image in $\Det(M)$ becomes a foliation chart of the canonical lift $\tilde{F}$ of $F$.  Thus $(\Det(M), \tilde{F})$ admits a finite covering of foliation charts which allows us to define the map $\Lambda\colon C^\infty_c(\Det(G)) \rightarrow H^*(\Det(M))$ in the same way as in the compact case.  

 This correspondence $\Lambda$ has two important special cases.  The first case concerns a holonomy invariant transverse measure (if any) of $(M; F)$.  In that case, the associated trace, which is a $0$-cyclic cocycle on $C_c(M; F)$, is mapped to the Poncar\'{e} dual $q$-form of the Ruelle--Sullivan $p$-current.  The second is that the transverse fundamental class is mapped to the class of unit in $H^0(M)$.
 
 Here, we are interested in the image of the $(q+1)$-cyclic cocycle $i_D \phi$ when $(M; F)$ admits a projectively invariant transverse measure.  From the consideration of parity, it follows that $\Lambda(i_D\phi)$ is a class odd degree.  The following theorem identifies this cohomology class.

\begin{thm}
  \label{thm:calcul-id-phi-as-grad-1-form}
  Suppose that there is a covering $\pi\colon M' \rightarrow M$ and a holonomy invariant transverse measure $\omega$ of $(M', \pi^* F)$ which is projectively invariant under the deck transformations of $\pi$.  Let $\phi$ be a transverse fundamental cocycle on $C^\infty_c(G_T)$, and $D$ be the derivation of the modular automorphism group associated with $\omega$.  Then $\Lambda(i_D \phi)$ agrees with $[\log \delta^{(\omega)}]$ in $H^1(M)$.
\end{thm}

In order to prove the above we take a well-behaved open sets of $M$ with
respect to the given foliation.

\begin{dfn}
  An open set $U$ in $M$ is said to be $F$-\textit{straight} when we have a
  complete transversal $T$ in $U$ with trivial holonomy for $F|_U$.
\end{dfn}

\begin{rmk}
  Since any foliation chart is $F$-straight, $M$ admits an covering by
  $F$-straight open sets.
\end{rmk}

Let $U$ be an $F$-straight open set of $M$.  Let $\tilde{U}$ be the preimage of $U$ with respect to the projection map $\tilde{\pi} \colon \Det(M) \rightarrow M$. Then $\tilde{U}$ is a $\tilde{\pi}^*F$-straight open set of $\Det(M)$.  Let $f$ be a section of the covering map.  Then the pullback of $\omega$ by $f$ determines a trivialization $\tilde{U} \simeq \R \times U$.  We let $t$ denote the coordinate on the first component of the right hand side.  Let $f^0$, $f^1$ be two sections $U \rightarrow M'$, and $c_{f^0 / f^1}$ be the scalar such that $c_{f^0/f^1} \omega = g^* \omega$ where $g$ is the deck transformation satisfying $g \circ f^0 = f^1$.  Then the $t$-coordinates determined respectively by $f^0$ and $f^1$ differ by a translation by the constant $\log (c_{f^0 / f^1})$.

\begin{proof}[Proof of Theorem~\ref{thm:calcul-id-phi-as-grad-1-form}]

We choose a covering $(U_i)_{i \in I}$ of $M$ by $F$-straight open sets, and let $(\tilde{U}_i)_{i \in I}$ the corresponding covering of $\Det(M)$ by the inverse images of the $U_i$.  For each $i$, we choose a section $f^i\colon U_i \rightarrow M'$ and let $t_i$ denote the corresponding $t$-coordinate on $\tilde{U}_i$.  The coordinate transforms among the $\tilde{U}_i$ are given by the ones for $U_i$, and the between the $t_i$-coordinates induced by the holonomy transformation.  We see that the $1$-forms $d t_i$ on the sets $\tilde{U}_i$ are invariant under these coordinate transforms.  Hence we obtain a global $1$-form on $\Det(M)$, denoted by $d t$.

Since the projection map $\Det(M) \rightarrow M$ has the contractible fiber $\R_+^\times$, it induces an isomorphism of the cohomology.  It is easy to see that the \v{C}ech $1$-cocycle corresponding to $d t$ is precisely the one used to define $[\log \delta^{(\omega)}]$.  We claim that $d t$ is the image of $\hat{\phi}$ under $\Lambda \colon \HP^*(C^\infty_c(\Det(G))) \rightarrow H^{* + q + 1}(\Det(M))$.

First, the restriction of $\hat{\phi}$ to $C^\infty_c(\tilde{U}_i ; \tilde{\pi}^*F)$ corresponds to the holonomy invariant transverse $q$-current on $\tilde{U}_i$ defined by
\[
\int_{\R} d t_i \int_{\ensemble{t} \times T_i} f^0 d f^1 \cdots d f^q.
\]
Integration by parts shows that this current is the boundary of the invariant transverse $(q+1)$-current
\[
\psi_i(f^0 d f^1 \cdots d f^{q+1}) = \int_{\tilde{U}_i} t_i f^0 d f^1 \cdots d f^{q + 1}.
\]
The \v{C}ech coboundary of $(\psi_i)_{i \in I}$ is given by $((t_i - t_j) \tilde{\phi}_{i j})_{i, j \in I}$, where $\tilde{\phi}_{i j}$ is the transverse fundamental current for $\tilde{\pi}^*F$.  This shows that $\Lambda(\hat{\phi})$ is indeed equal to the class of $d t$.

By Lemma~\ref{lem:topol-desc-connes-thom}, we obtain that $i_D \phi$ and $\hat{\phi}$ induce the same map on the geometric $K$-groups.  From this it already follows that $\Lambda(i_D \psi)$ and $[\log \delta^{(\omega)}]$ pairs the same way against the homology classes of the form $\Todd(\tau \otimes \C)^{-1}\ch(\sigma_D) \cap [M]$ for a longitudinal elliptic operator $D$.  To see that they agree as elements of $H^*(M)$, we may use the fact that the construction of $\hat{\phi}|_{C^\infty_c(\tilde{U}; \tilde{\pi}^*F)}$ and the natural isomorphisms $\HP^*(C^\infty_c(\tilde{U}; \tilde{\pi}^*F)) \rightarrow \HP^{*+1}(C^\infty_c(U; F))$ for the open sets $U$ of $M$ are compatible with the localization maps used in the definition of $\Lambda$.
\end{proof}

\begin{rmk}
The setting of this section is already considered by H. Moriyoshi~\cite{MR2285046}, who raised a question regarding the relationship between the structure of $W(M; F)$ and a certain geometrically defined subset of $\R$, called the \textit{$K$-set} of $(M; F)$.

Let $H_1^F$ be the intersection of $H_1(M; \Z)$ and the joint kernel of the closed $F$-basic $1$-forms.  Intuitively, $H_1^F$ can be thought of as the span of the $1$-cycles in the leaves of $F$.  Then, the $K$-set of $(M; F)$ is defined as the values which $[\log \delta^{(\omega)}]$ takes on $H_1^F$.

Suppose that the von Neumann algebra $W(M; F)$ is of type III$_\lambda$.  Theorem~\ref{thm:calcul-id-phi-as-grad-1-form}, combined with Theorem~\ref{thm:iii-lambda-proj-inv-per-calc-trans-fund}, and the compatibility of $\Lambda$ with $\Phi_K$ shows that the pairing of $[\log \delta^{(\omega)}]$ with the twisted Chern character of the cycles in $K^{\topol}_*(M; F)$ is contained in $\Q \log(\lambda)$.  If the Baum--Connes conjecture for $(M; F)$ holds, the number $\log(\lambda)$ itself can be realized by some cycle.  The set $\pairing{i_D}{\phi, K^{\topol}_*(M; F)}$ will contain the $K$-set if the twisted Chern character of an element $y \in K^{\topol}_*(M; F)$ satisfying $\partial(\mu(y)) = [G_T^{(u)}]^*$ (see the proof of Theorem~\ref{thm:iii-lambda-proj-inv-per-calc-trans-fund}) lies in $H^F_1$, which happens in many examples.
\end{rmk}

\paragraph{Acknowledgements}  This paper was a part of author's master thesis
submitted to the University of Tokyo.  He would like to thank Y. Kawahigashi for
his continuing support throughout the research.  He is also grateful to
G. Skandalis, C. Oikonomides, H. Moriyoshi, T. Katsura, and T. Tsubo\^{\i} for
many fruitful exchanges.  Lastly but not least, he thanks the referee for his valuable comments which were crucial to fix many ambiguities in the original manuscript.


\end{document}